\documentclass[11pt, reqno]{amsart}
\usepackage{amsmath}
\usepackage{amssymb}
\usepackage{amsfonts}
\usepackage[usenames]{color}
\usepackage{version}

 \newtheorem{theorem}{Theorem}[section]
 \newtheorem{corollary}[theorem]{Corollary}
 \newtheorem{lemma}[theorem]{Lemma}
 \newtheorem{proposition}[theorem]{Proposition}

\newtheorem{definition}[theorem]{Definition}
\newtheorem{example}[theorem]{Example}

\newtheorem{fact*}{Fact}

\DeclareMathOperator{\RE}{Re}
\DeclareMathOperator{\IM}{Im}

\DeclareMathOperator{\dist}{dist}

\DeclareMathOperator{\spa}{span}
\DeclareMathOperator\re{\mathrm {Re~}}
\DeclareMathOperator\im{\mathrm {Im~}}
\newcommand\dd{\mathrm d}

\newcommand{\s}{\mathcal{S}_2}
\newcommand{\sd}{\mathcal{S}_d}

\newcommand\Ha{\mathbb{H}}
\renewcommand{\L}{\mathcal{L}}
\newcommand{\M}{\mathcal{M}}
\newcommand{\N}{\mathcal{N}}
\newcommand{\T}{\mathbb{T}}

\newcommand{\D}{\mathbb{D}}
\newcommand{\C}{\mathbb{C}}

\newcommand{\R}{\mathbb{R}}

\newcommand{\cc}[1]{\overline{#1}}
\newcommand{\abs}[1]{\left\vert#1\right\vert}
\newcommand{\set}[1]{\left\{#1\right\}}

\newcommand{\norm}[1]{\left\Vert#1\right\Vert}

\newcommand{\ran}[1]{\operatorname{ran}#1}
\newcommand{\nt}{\stackrel{\mathrm {nt}}{\to}}

\newcommand{\ip}[2]{\left\langle #1, #2 \right\rangle}
\newcommand{\ad}{^\ast}

\newcommand\Pick{\mathcal P}
\newcommand{\p}{\mathcal{P}_2}
\newcommand{\HH}{\Ha\times\Ha}
\newcommand\Htau{\mathbb{H}(\tau)}
\newcommand{\vp}{\varphi}
\newcommand{\ph}{\varphi}
\newcommand\al{\alpha}
\newcommand\ga{\gamma}
\newcommand\de{\delta}

\newcommand\la{\lambda}
\renewcommand\ln{{\la_n}}

\newcommand\si{\sigma}
\newcommand\beq{\begin{equation}}
\newcommand\ds{\displaystyle}
\newcommand\eeq{\end{equation}}
\newcommand\df{\stackrel{\rm def}{=}}
\newcommand\ii{\mathrm i}
\newcommand{\vectwo}[2]
{
   \begin{pmatrix} #1 \\ #2 \end{pmatrix}
}
\newcommand\half{\tfrac 12}
\newcommand\blue{\color{blue}}
\newcommand\black{\color{black}}
\newcommand\red{\color{red}}
\newcommand\nn{\nonumber}
\newcommand\nin\noindent
\numberwithin{equation}{section}

\begin{document}

\title[Boundary behavior via generalized models]{Boundary behavior of analytic functions of two variables via generalized models}

\author{J. Agler, R. Tully-Doyle and N. J. Young}

\thanks{The first author was partially supported by National Science Foundation Grant on  Extending Hilbert Space Operators DMS 1068830. The third author was partially supported by the UK Engineering and Physical Sciences Research Council grant EP/J004545/1.}

\date{26th March 2012}
\begin{abstract}
We describe a generalization of the notion of a Hilbert space model of a function in the Schur class of the bidisc.  This generalization is well adapted to the investigation of boundary behavior at a mild singularity of the function on the $2$-torus.  We prove the existence of a generalized model with certain properties corresponding to such a singularity and use this result to solve two function-theoretic problems.  The first of these is to characterise the directional derivatives of a function in the Schur class at a singular point on the torus for which the Carath\'eodory condition holds.  The second is to obtain a representation theorem for functions in the two-variable Pick class analogous to the refined Nevanlinna representation of functions in the one-variable Pick class.
\end{abstract}

\maketitle
\section{Introduction}\label{intro}
In this paper we solve two problems about analytic functions of two variables using a variant of the notion of a Hilbert space model of a function.  One problem concerns the generalization to two variables of a classical representation theorem of Nevanlinna, while the other is quite unlike any question that arises for functions of a single variable.  Both relate to behavior of functions at boundary points of their domains.

The first problem is: what directional derivatives are possible for a function in the two-variable Schur class $\s$ at a singular point on the $2$-torus $\T^2$?  To clarify this question let us consider the rational function
\beq\label{favourite}
\ph(\la) = \frac{ \half \la_1+\half \la_2 -\la_1\la_2}{1-\half \la_1-\half\la_2}, \quad \la\in\D^2,
\eeq
where $\D$ denotes the open unit disc.  This function belongs to $\s$ (that is, it is analytic and bounded by $1$ in modulus on $\D^2$).  It has a singularity at the point $\chi=(1,1)\in\T^2$, in that $\ph$ does not extend analytically (or even continuously) to $\chi$.  Nevertheless $\ph$ has nontangential limit $1$ at $\chi$, and so we may define $\ph(\chi)$ to be $1$.  Despite the fact that $\ph$ is discontinuous at $\chi$, the directional derivative $D_{-\de}\ph(\chi)$ exists for every direction $-\de$ pointing into the bidisc at $\chi$, and
\begin{align}\label{DirDer}
D_{-\de}\ph(\chi) &=-\frac{2\de_1\de_2}{\de_1+\de_2}  \nn\\
	&= \ph(\chi)\de_2 h(\de_2/\de_1)
\end{align}
where $h(z)=-2/(1+z)$.  

Remarkably enough, a similar statement holds in great generality \cite[Theorem 2.10]{AMY10}.  If $\ph\in\s$ has a singularity at $\chi$ and $\ph$ satisfies a weak regularity condition at $\chi$ (the Carath\'eodory condition, explained in Section \ref{carapoints}) then $D_{-\de}\ph(\chi)$ exists for all relevant directions $\de$, and  furthermore there exists an analytic function $h$ on the upper halfplane
\[
\Pi=\{z : \im z > 0\}
\]
such that both $h(z)$ and $-zh(z)$ have non-negative imaginary part and the directional derivative $D_{-\de}\ph(\chi)$ is given by equation \eqref{DirDer}.  We call $h$ the {\em slope function} for $\ph$ at $\chi$.  The problem, then, is to find necessary and sufficient conditions for a function $h$ on $\Pi$ to be the slope function of some member of $\s$.   It transpires that the stated necessary conditions on $h$ are also sufficient for $h$ to be a slope function (Theorem \ref{SuffSlope} below).

The second problem is to generalize to two variables a theorem of Nevanlinna which plays an important role in one proof of the spectral theorem for self-adjoint operators.  Nevanlinna's theorem gives an integral representation formula for the functions in the Pick class $\Pick$ (the analytic functions on $\Pi$ having non-negative imaginary part) that satisfy a growth condition on the imaginary axis; it states that such functions are the Cauchy transforms of the finite positive measures on the real line $\R$.  Nevanlinna's growth condition can be regarded as a regularity condition at the point $\infty$ on the boundary of $\Pi$.  We obtain an analogous representation for functions in the two-variable Pick class that satisfy a suitable regularity condition at $\infty$, but rather than an integral formula we get an expression involving the two-variable resolvent of a densely defined self-adjoint operator on a Hilbert space (Theorem \ref{NevanlinnaRep}).

To solve these two problems we modify the notion of model so as to focus on the behavior of a function  $\ph\in\s$ near a boundary point at which $\ph$ satisfies Carath\'eodory's condition (see Definition \ref{defCara} below).
A {\em model} of an analytic function $\ph$ on the polydisc $\D^d$ is a pair $(\M,u)$ where $\M$ is a separable Hilbert space with an orthogonal decomposition $\M=\M_1\oplus\dots\oplus \M_d$ and $u:\D^d \to \C$ is an analytic map such that, for all $\la, \mu \in \D^d$,
\beq\label{genModEq}
1-\overline{\ph(\mu)} \ph(\la) = \ip{(1-I(\mu)^*I(\la))u_\la}{u_\mu},
\eeq
where $I(\la) = \la_1P_1+\dots +\la_d P_d$ and $P_j$ is the orthogonal projection on $\M_j$.  This notion is particularly effective in the case $d=2$, since every function in $\s$ has a model \cite{Ag1}.   We used models, and their accompanying realizations, in \cite{AMY10} to prove a Carath\'eodory theorem for functions in $\s$, but for our present purpose it is too restrictive to require that $I(\la)$ in equation \eqref{genModEq} be linear in $\la$.  By allowing $I(.)$ to be a general operator-valued inner function on $\D^2$ we acquire greater flexibility.   In Theorem \ref{mainExistence} we prove the existence of a model of $\ph\in\s$, of this more general type, with special properties relative to a boundary point at which $\ph$ satisfies Carath\'eodory's condition.
Such generalized models then provide the main tool for the solution of our two problems.

The paper is organized as follows.  In Section \ref{carapoints} we recall some definitions and discuss the Carath\'eodory condition.  In Section \ref{GenModel} we define generalized models and prove the main existence theorem for them.  In Section \ref{DirDerivs} we use generalized models to give an alternative proof of the existence of directional derivatives and slope functions.   In Section \ref{represent} we derive an integral representation formula for functions $h$ on $\Pi$ such that both $h$ and $-zh$ belong to the Pick class, and in Section \ref{converse} we use this integral representation to construct a function in $\s$ having slope function $h$ at a point on the torus.  In Section \ref{NevanRep} we prove a two-variable analog of Nevanlinna's representation theorem for functions in the Pick class subject to a growth condition on the imaginary axis.

\section{Carapoints}\label{carapoints}
C. Carath\'eodory in \cite{cara} proved that if a function $\ph$ in the one-variable Schur class satisfies 
\beq\label{caraCond}
\liminf_{\la\to\tau} \frac{1-|\ph(\la)|}{1-|\la|} < \infty
\eeq
for some $\tau\in\T$ then not only does $\ph$ have a nontangential limit at $\tau$, but it also has an angular derivative $\ph'(\tau)$ at $\tau$, and $\ph'(\la)\to\ph'(\tau)$ as $\la$ tends nontangentially to $\tau$ in $\D$.
Here nontangential limits are defined as follows.  For any domain $U$ and for $\tau$ in the topological boundary $\partial U$ of $U$ we say that a set $S\subset U$ {\em approaches $\tau$ nontangentially} if $\tau\in S^-$, the closure of $S$, and
\[
\left\{ \frac{\|\la-\tau\|}{\dist(\la,\partial U)} : \la\in S\right\} \mbox{ is bounded}.
\]
We say that a function $\ph$ on $U$ {\em has nontangential limit $\ell$ at $\tau$}, in symbols
\[
\lim_{\la\nt\tau} \ph(\la)=\ell,
\]
if
\[
\underset{\la\in S}{\lim_{\la\to\tau}}  \, \ph(\la) = \ell
\]
for every set $S\subset U$ that approaches $\tau$ nontangentially.

Carath\'eodory's result has been generalized by several authors, notably by K. W\l odarczyk \cite{woody}, W. Rudin \cite{rudin}, F. Jafari \cite{jaf93},  M. Abate \cite{abate} and two of us with J. E. McCarthy \cite{AMY10}.  Carath\'eodory's condition \eqref{caraCond} generalizes naturally to holomorphic maps $\ph:U\to V$ for any pair of bounded domains $U, V$ in complex Euclidean spaces of finite dimensions.  For any $\tau$ in $\partial U$ we say that $\ph$ {\em satisfies the Carath\'eodory condition at} $\tau$, or that $\tau$ is a {\em carapoint} for  $\ph$, if 
\beq\label{defCara}
\underset{\la\in U}{\liminf_{\la\to \tau}} \,  \frac{\dist(\ph(\la), \partial V)}{\dist(\la,\partial U)} < \infty.
\eeq
  In particular, when $U=\D^d, \, V=\D$ and $\ph\in\sd$, $\tau$ is a carapoint for $\ph$ if
\[
\liminf_{\la\to\tau}\frac{1-|\ph(\la)|}{1-\|\la\|_\infty} < \infty.
\]
Likewise, if $\ph$ is a contractive operator-valued analytic function on $\D^d$, $\tau\in\T^d$ is a carapoint for $\ph$ if
\[
\liminf_{\la\to\tau}\frac{1-\|\ph(\la)\|}{1-\|\la\|_\infty} < \infty.
\]
Of course any point in $\T^d$ at which $\ph$ is analytic is a carapoint for $\ph$, but we are concerned here with {\em singular} carapoints. We say that an analytic function $\ph$ on a domain $U$ is {\em singular} at a point $\tau\in\partial U$ if there is no neighborhood $W$ of $\tau$ such that $\ph$ extends to an analytic function on $U\cup W$.

In Section \ref{NevanRep} we shall also define carapoints at infinity for certain unbounded domains $U$ and $V$.

Not all the conclusions of Carath\'eodory's Theorem hold even for $\s$:  for $\ph$ of the example \eqref{favourite} of Section \ref{intro}, $\chi$ is a carapoint, but since $D_{-\de}\ph(\chi)$ is not linear in $\de$, it is not the case that $\ph$ has an angular gradient at $\chi$.  Indeed, the interest in the first of our two problems is precisely in carapoints at which there is no angular gradient.  However it is true for all the cases considered in this paper that if $\tau\in\partial U$ is a carapoint for $\ph:U\to V$ then $\ph$ has a nontangential limit at $\tau$ \cite{woody}.  This limit will be denoted by $\ph(\tau)$; it is obvious that $\ph(\tau)\in\partial V$.

Here is some more terminology and notation.  We denote by $\Ha$ the right halfplane $\{z\in\C: \re z> 0\}$.  
An operator-valued analytic function $I$ on $\D^d$ is said to be {\em inner} if $I(\la)$ is a unitary operator for almost all $\la\in\T^d$ with respect to Lebesgue measure.    The {\em Schur class of the polydisc $\D^d$} is the set of analytic functions from $\D^d$ to the closed unit disc $\D^-$ and is denoted by $\mathcal{S}_d$.

\section{Generalized models of Schur-class functions} \label{GenModel}

In the definition of a model of a function $\ph:\D^d\to\C$  (see equation \eqref{genModEq} above), the co-ordinate functions have a privileged position through the definition of $I(.)$ as linear in the co-ordinates.  One consequence is that any singular behavior of $\ph$ at a boundary point must be reflected in singular behavior of $u$, rather than $I(.)$, near that point.  A simple relaxation of the definition of model enables us to concentrate information about singular behavior in the inner function $I(.)$ instead, and this proves helpful for the two problems we study here.

\begin{definition}\label{genmodel}
Let $\vp:\D^d\to\C$ be analytic. The triple $(\mathcal M, u, I)$ is a {\em generalized model of $\vp$} if 
\begin{enumerate}
 \item $\M$ is a separable Hilbert space, 
 \item $u: \D^d \to \M$ is analytic, and
 \item $I$ is a contractive analytic  $\mathcal{L}(\M)$-valued function on $\D^d$
\end{enumerate}
  such that the equation
\begin{equation}\label{modeleqn}
1- \cc{\vp(\mu)}\vp(\lambda) = \ip{(1 - I(\mu)^* I(\lambda))u_\lambda}{u_\mu}
\end{equation}
holds for all $\lambda, \mu \in \D^d$.

The generalized model $(\mathcal M, u, I)$ is {\em inner} if $I(.)$ is inner.
\end{definition}

Clearly, in the case that $I(\lambda) = \lambda_1 P_1 + \dots +\lambda_d P_d$, we recapture the notion of model in the previous sense.

A well-known lurking isometry argument proceeds from a model $(\M,u)$ of a function $\ph\in\sd$ to a realization of $\ph$ \cite{Ag1}.  The identical argument applied to a generalized model $(\M, u, I)$ produces a generalized notion of realization.

\begin{theorem} \label{genRealiz}
 If $(\L, u, I)$ is a generalized model of $\vp \in \sd$ then there exist a Hilbert space $\M$ containing $\L$, a scalar $a \in \C$, vectors $\beta, \gamma \in \M$ and a linear operator $D:\M \to \M$ such that the operator
\beq \label{defL}
L =  \begin{bmatrix}
  a & 1 \otimes \beta \\
  \gamma \otimes 1 & D
 \end{bmatrix}
\eeq
is unitary on $\C \oplus \M$ and, for all $\lambda \in \D^d$, 
\beq\label{Leq}
 L \vectwo{1}{I(\lambda)u_\lambda} = \vectwo{\vp(\lambda)}{u_\lambda},
\eeq
and consequently, for all $\la\in\D^d$,
\beq\label{formphi}
\ph(\la) = a+ \ip{I(\la)(1-DI(\la))^{-1}\ga}{\beta}.
\eeq
\end{theorem}
\begin{proof}
By equation \eqref{modeleqn}, for all $\la,\mu\in\D^d$,
\[
1+\ip{I(\la)u_\la}{I(\mu)u_\mu} = \overline{\ph(\mu)}\ph(\la) + \ip{u_\la}{u_\mu}.
\]
We may interpret this equation as an equality between the gramians of two families of vectors in $\C\oplus\L$.
Accordingly we may define an isometric operator
\[
L_0: \spa \left\{ \vectwo{1}{I(\lambda)u_\lambda} : \la\in\D^d\right\}\to \spa\left\{ \vectwo{\vp(\lambda)}{u_\lambda}: \la\in\D^d \right\}
\]
by equation \eqref{Leq}.  If necessary we may enlarge $\C\oplus\L$ to a space $\C\oplus\M$ in which the domain and range of $L_0$ have equal codimension, and then we may extend $L_0$ to a unitary operator $L$ on $\C\oplus\M$.
\end{proof}
The ordered 4-tuple $(a, \beta, \gamma, D)$, as  in equation \eqref{defL}, will be called a {\em realization} of the (generalized) model $(\L, u, I)$ of $\ph$ if $L$ is a contraction and equation \eqref{Leq} holds.  It will be called a {\em unitary realization} if in addition $L$ is unitary on $\C\oplus \L$. 

Realizations provide an effective tool for the study of boundary behavior.  Here is a preliminary observation.
\begin{lemma}\label{Nperp}
 Suppose that $\vp \in \sd$ has a model $(\M, u)$ with realization $(a, \beta, \gamma, D)$. Let $\tau\in\T^d$ be a carapoint for $\ph$ and let $\mathcal N = \ker(1 - D\tau)$.  Then
\beq\label{gabeta} 
\gamma \in \ran (1-D\tau) \subset\N^\perp \mbox{  and } \tau^*\beta \in \N^\perp.
\eeq 
\end{lemma}
\begin{proof}
 First we show that $\tau^*\beta \in \mathcal N^\perp$. Let $L$ be given by equation \eqref{defL}. Choose any $x \in \mathcal N$.   Then $x=D\tau x$ and so
\[
L\vectwo{0}{\tau x} =  \begin{bmatrix}
  a & 1 \otimes \beta \\
  \gamma \otimes 1 & D
 \end{bmatrix} \vectwo{0}{\tau x} = \vectwo{\ip{\tau x}{ \beta}}{D\tau x} = \vectwo{\ip{x}{\tau^*\beta}}{x}.
\]
 Since $L$ is a contraction and $\tau$ is an isometry,
\[
\norm{ \vectwo{\ip{x}{\tau^*\beta}}{x}} = \norm{L\vectwo{0}{\tau x}} \leq  \norm{\tau x} = \norm{x}, 
\]
and so $\ip{x}{\tau^*\beta} = 0$. Since $x\in \mathcal N$ is arbitrary, $\tau^*\beta \in \N^\perp$.

Proposition 5.17 of \cite{AMY10} asserts that $\tau$ is a carapoint for $\ph$ if and only if $\ga\in \ran(1 - D\tau)$. Now since $D\tau$ is a contraction, every eigenvector of $D\tau$ corresponding to an eigenvalue $\la$ of unit modulus is also an eigenvector of $(D\tau)^*$ with eigenvalue $\bar\la$.  Hence $\ker(1-D\tau)=\ker(1-\tau^*D^*)$, and we have
\[
 \gamma \in \ran(1-D\tau) \subset \ker(1-\tau^*D\ad)^\perp = \ker(1-D\tau)^\perp = \N^\perp.
\]
\end{proof}

We are interested in the behavior of models at carapoints of $\ph\in\sd$.  Here are two relevant notions.
\begin{definition} \label{defBCpts}
Let $(\M, u, I)$ be a generalized model of a function $\ph\in\sd$.  A point $\tau\in\partial\D^d$ is a {\em $B$-point} of the model if $u$ is bounded on every subset of $\D^d$ that approaches $\tau$ nontangentially.  The point $\tau$ is a {\em $C$-point} of the model if, for every subset $S$ of $\D^d$ that approaches $\tau$ nontangentially, $u$ extends continuously to $S \cup \set{\tau}$ (with respect to the norm topology of $\M$).
\end{definition}

As is well known, not all functions in $\sd$ have models when $d\geq 3$.  For the rest of the paper we restrict attention to the case $d=2$; in this case it {\em is} true that every function in the Schur class has a model \cite{Ag1}.

Our next task is to show that if a function  $\vp\in\s$ has a singularity at a $B$-point $\tau$, then we can construct a generalized model of $\ph$ in which the singularity of $\ph$ is encoded in an $I(\lambda)$ that is singular at $\tau$, in such a way that the model has a $C$-point at $\tau$.  The device that leads to this conclusion is to write vectors in and operators on $\M$ in terms of the orthogonal decomposition $\M=\N\oplus\N^\perp$ where $\N= \ker(1-D\tau)$ and $D$ comes from a realization of $(\M,u)$.  The following observation is straightforward.
\begin{lemma}\label{decompP1}
Let $\N$ be a subspace of $\M$ and let $P_1$ be a Hermitian projection on $\M$.  With respect to the decomposition $\N\oplus\N^\perp$ the operator $P_1$ has operator matrix
\beq\label{formP1}
P_1 = \begin{bmatrix}
              X & B \\ B\ad & Y
             \end{bmatrix}
\eeq
for some operators $X, Y, B$, where
\begin{enumerate}
 \item  $0 \leq X, Y\leq 1$
 \item $BB\ad = X(1-X),  \quad B\ad B = Y(1 - Y)$
 \item $BY = (1-X)B, \quad B(1 - Y) = XB$
 \item $B\ad X = (1 - Y)B\ad, \quad B\ad(1 - X) = YB\ad$.
\end{enumerate}
\end{lemma}
We now construct a generalized model corresponding to a carapoint of $\ph\in\s$.
\begin{theorem}\label{mainExistence}
Let $\tau\in\T^2$ be a carapoint for $\ph\in\s$.   There exists an inner generalized model $(\M,u,I)$ of $\ph$ such that 
\begin{enumerate}
\item $\tau$ is a $C$-point for $(\M,u,I)$, 
\item $I$ is analytic at every point $\la \in\T^2 $ such that $\la_1\neq \tau_1$ and $\la_2\neq \tau_2$, and 
\item $\tau$ is a carapoint for $I$ and $I(\tau)=1_\M$.
\end{enumerate}
Furthermore, we may express $I$ in the form
\begin{equation}\label{defI}
  I(\lambda) = \frac{\bar\tau_1\la_1 Y + \bar\tau_2\la_2(1-Y)- \bar\tau_1\bar\tau_2 \la_1\la_2}{1-\bar\tau_1\la_1(1-Y)-\bar\tau_2\la_2 Y}
 \end{equation}
for some positive contraction $Y$ on $\M$.
\end{theorem}
\begin{proof}
Choose any model $(\L,v)$ of $\ph$ and any realization $(a,\beta_0,\gamma,D)$ of $(\L,v)$.  By definition, $\L$ comes with an orthogonal decomposition $\L=\L_1 \oplus \L_2$: let $P_1$ be the orthogonal projection on $\L_1$.  Since $\tau$ is a carapoint for $\ph$ we may apply Lemma \ref{Nperp} to deduce that $\ga\in\ran(1-D\tau)$ and $\tau^*\beta_0, \ga \in\ker(1-D\tau)^\perp$.

Consider first the case that $\ker(1-D\tau) = \{0\}$.  This relation  implies that there is a {\em unique} vector $v_\tau \in\L$ such that $(1-D\tau)v_\tau = \ga$.   Let $(\ln)$ be any sequence in $\D^2$ that converges nontangentially to $\tau$.  We claim that $v_\ln \to v_\tau$.  Suppose not: then since $(v_\ln)$ is bounded, by \cite[Corollary 5.7]{AMY10}, we can assume on passing to a subsequence that $(v_\ln)$ tends weakly to a limit $x\in\L$ different from $v_\tau$.  By \cite[Proposition 5.8]{AMY10} it follows that $v_\ln \to x$ in norm.  Take limits in the equation
\[
(1-D\ln)v_\ln = \ga
\]
to deduce that $(1-D\tau) x =\ga$.  Since $x\neq v_\tau$, this contradicts the fact that $(1-D\tau)^{-1}\ga=\{v_\tau\}$.  Hence $v_\ln \to v_\tau$.  In other words $v$ extends continuously to $S\cup \{\tau\}$ for any set $S$ in $\D^2$ that tends nontangentially to $\tau$, which is to say that $\tau$ is a $C$-point for the model $(\L,v)$.  The conclusion of the theorem therefore holds if we simpy take $\M=\L, \, u=v$ and $I(\la)= \la_1P_1 + \la_2 P_2$.

Now consider the case that $\ker(1-D\tau)\neq \{0\}$.
 Let $\N= \ker (1-D\tau)$. 
With respect to the decomposition $\L = \N\oplus \N^\perp$ we may write
\beq\label{defQ}
 D\tau = \begin{bmatrix} 1 & 0 \\ 0 & Q \end{bmatrix}
\eeq
and $v_\lambda = \vectwo{w_\lambda}{u_\lambda}$.   Note that $\ker (1-Q) = \{0\}$.

Let us express $\la $, acting as an operator on $\L$ by
\[
 \lambda = \lambda_1 P_{1} \oplus \lambda_2 (1 - P_{1}),
\]
as an operator matrix with respect to the decomposition $\L=\N\oplus\N^\perp$, as in Lemma \ref{decompP1}:
\beq\label{formla}
 \lambda =\la_1P_1 + \la_2(1-P_1) =  \begin{bmatrix}
            \lambda_1 X + \lambda_2 (1-X) & (\lambda_1 - \lambda_2)B \\ 
            (\lambda_1 - \lambda_2)B\ad & \lambda_1 Y + \lambda_2 (1-Y)
           \end{bmatrix}
\eeq
where $X,\ Y$ are the compressions of $P_1$ to $\N, \N^\perp$ respectively, so that $0 \leq X,Y \leq 1$.
Thus
\[
 1 - D\lambda =1-D\tau \tau^*\la= \begin{bmatrix}
                 1 - (\lambda'_1 X + \lambda'_2 (1-X)) & -(\lambda'_1 - \lambda'_2)B \\
                 - (\lambda'_1 - \lambda'_2)QB\ad & 1 - Q(\lambda'_1 Y + \lambda'_2 (1 - Y))
                \end{bmatrix}
\]
where $\la'_1=\bar\tau_1\la_1, \ \la'_2=\bar\tau_2\la_2$.
Since $(1 - D\lambda)v_\lambda = \gamma$,
\[
 \vectwo{0}{\gamma} =(1 - D\lambda)v_\lambda = \begin{bmatrix}
                 1 - (\lambda'_1 X + \lambda'_2 (1-X)) & -(\lambda'_1 - \lambda'_2)B \\
                 - (\lambda'_1 - \lambda'_2)QB\ad & 1 - Q(\lambda'_1 Y + \lambda'_2 (1 - Y))
                \end{bmatrix}\vectwo{w_\lambda}{u_\lambda},
\]
from which we have the equations
\begin{align}
 (1 - \lambda'_1 X - \lambda'_2 (1-X))w_\lambda - (\lambda'_1 - \lambda'_2)Bu_\lambda = 0 \label{eq:a1}\\
 -(\lambda'_1 - \lambda'_2)QB\ad w_\lambda + (1 - Q(\lambda'_1 Y + \lambda'_2 (1-Y)))u_\lambda = \gamma \label{eq:a2}.
\end{align}

By Lemma \ref{decompP1} and equation \eqref{eq:a1}  we have
\begin{align} \label{propvla}
 0 &=B^*\left((1 - \la'_1 X - \la'_2 (1-X))w_\lambda - (\la'_1 - \la'_2)Bu_\lambda\right) \nn \\
 &= (B\ad - \la'_1 B\ad X - \la'_2 B\ad (1 - X)) w_\lambda - (\la'_1 - \la'_2)B\ad B u_\lambda  \nn\\
 &= (B\ad - \la'_1 (1-Y)B\ad - \la'_2 Y B\ad ) w_\lambda - (\la'_1 - \la'_2)Y(1- Y) u_\lambda \nn\\
 &= (1 - \la'_1 (1-Y) - \la'_2 Y)B\ad w_\lambda - (\la'_1 - \la'_2) Y(1 - Y)u_\lambda.
\end{align}
Since $0\leq Y\leq 1$ it is clear from the spectral mapping theorem that
\[
1\notin \si(\la'_1(1-Y)+\la'_2 Y)
\]
for all $\la\in\D^2$, and thus equation \eqref{propvla} tells us that
\begin{equation}\label{Bstarv}
 B\ad w_\lambda = \frac{(\la'_1 - \la'_2)Y(1 - Y)}{1 - \la'_1(1 - Y) - \la'_2 Y}u_\lambda.
\end{equation}

Substituting the relation \eqref{Bstarv} into \eqref{eq:a2} we obtain
\begin{align}
 \ga &=-(\la'_1 - \la'_2)QB\ad w_\lambda + (1 - Q(\la'_1 Y +\la'_2 (1-Y))u_\lambda \notag \\
 &= -(\la'_1 - \la'_2)Q \frac{(\la'_1 - \la'_2)Y(1 - Y)}{1 - \la'_1(1 - Y) - \la'_2 Y}u_\lambda + (1 - Q(\la'_1 Y +\la'_2 (1-Y))u_\lambda \notag \\
 &= \left[ 1 - Q\left(\frac{(\la'_1 - \la'_2)^2 Y(1 - Y)}{1 - \la'_1(1 - Y) - \la'_2 Y} + \la'_1 Y +\la'_2 (1-Y)\right)\right]u_\lambda \notag \\
 &= (1 - QI(\lambda))u_\lambda  \label{eq:b2}
\end{align}
where 
\begin{align}\label{gotI}
 I(\lambda) &= \frac{\la'_1 Y + \la'_2 (1-Y) - \la'_1 \la'_2}{1 - \la'_1(1-Y) - \la'_2 Y}  \\
   &= \frac{\bar\tau_1\la_1 Y + \bar\tau_2\la_2(1-Y)- \bar\tau_1\bar\tau_2 \la_1\la_2}{1-\bar\tau_1\la_1(1-Y)-\bar\tau_2\la_2 Y} \in \L(\M), \nn
\end{align}
which agrees with equation \eqref{defI}.  

Let  $\M=\N^\perp$: we claim that $(\M,u,I)$ is an inner generalized model of $\ph$ having the properties described in Theorem \ref{mainExistence}.

Firstly, it is clear from the formula \eqref{gotI} that $I$ is analytic on $\D^2$ and at every point $\la \in\T^2$ such that $1 \notin \si(\la'_1(1-Y) + \la'_2 Y)$. By the spectral mapping theorem and the fact that $0\leq Y \leq 1$, the spectrum of $\la'_1(1-Y) + \la'_2 Y$ is contained in the convex hull of the points $\la'_1, \ \la'_2$.  Hence $\si(\la'_1(1-Y) + \la'_2 Y)$ contains the point $1$ if and only if either $\la'_1=1$ and $0\in\si(Y)$ or $\la'_2=1$ and $1\in\si(Y)$.  Thus $I$ is analytic at points $\la\in\T^2$ for which $\la'_1\neq 1, \, \la'_2\neq 1$, that is, such that $\la_1\neq \tau_1, \, \la_2\neq \tau_2$.  The function $I$ therefore satisfies condition (2) of the theorem.

We must show that $I$ is an inner function.  Indeed, if $d(\la)$ denotes the denominator of  $I(\la)$ in equation \eqref{gotI}, we find that, for all $\la\in\Delta^2$ such that $1\notin \si(\la'_1(1-Y) + \la'_2 Y)$,
\begin{align*}
d(\la)^*(1-&I(\la)^*I(\la))d(\la) = \\
		&|1-\la'_1|^2(1-|\la'_2|^2) + 2\{\re (\la'_1 - \la'_2)-|\la'_1|^2 +|\la'_2|^2 +\re(\overline{\la'_1}\overline{\la'_2} (\la'_1-\la'_2)) \} Y.
\end{align*}
Hence $I(\la)^*I(\la)=1_\M$ for  all $\la\in\T^2$ such that $\la_1\neq \tau_1, \, \la_2\neq \tau_2$, and therefore for almost all $\la\in\T^2$ with respect to $2$-dimensional Lebesgue measure on $\T^2$.  Since $I(\la)$
 is clearly a normal operator for all such $\la$, it follows that $I$ is an inner $\L(\M)$-valued function.

Next we prove the model relation \eqref{modeleqn} for $(\M,u,I)$.   Let us calculate $\tau^*\lambda v_\lambda$  using equation \eqref{formla}:
\begin{align}\label{t*lvl}
\tau^*\la v_\la &= \tau^*\la\vectwo{w_\lambda}{u_\lambda} =\vectwo{(\la'_1 X+\la'_2(1-X))w_\la + (\la'_1-\la'_2)Bu_\la}{(\la'_1-\la'_2)B^*w_\la +(\la'_1Y+\la'_2(1-Y))u_\la}_{\N\oplus\N^\perp}.
\end{align}
By equations \eqref{Bstarv} and \eqref{t*lvl},
\begin{align*}
P_{\N^\perp} \tau^*\lambda v_\lambda 
		&= \left(\frac{(\la'_1 - \la'_2)^2 Y(1 - Y)}{1 - (\la'_1(1 - Y) + \la'_2 Y)} + \la'_1 Y +\la'_2 (1-Y)\right) u_\lambda \\
  &= I(\lambda)u_\lambda.
 \end{align*}
By equation \eqref{eq:a1},
\[
(\la'_1 X+\la'_2(1-X))w_\la = w_\la - (\la'_1-\la'_2)Bu_\la,
\]
which, in combination with equation \eqref{t*lvl}, yields the relation 
\[
P_\N \tau^*\la v_\la =(\la'_1 X+\la'_2(1-X))w_\la + (\la'_1-\la'_2)Bu_\la= w_\la
\]
and therefore
\[
\tau^*\la v_\la = \vectwo{w_\la}{I(\la)u_\la}_{\N\oplus\N^\perp}.
\]

Hence
\begin{align*}
1-\overline{\ph(\mu)}\ph(\la) &= \ip{(1-\mu^*\la)v_\la}{v_\mu}_\L \\
		&= \ip{v_\la}{v_\mu}_\L - \ip{\la v_\la}{\mu v_\mu}_\L \\
		&= \ip{w_\la}{w_\mu}_\N + \ip{u_\la}{u_\mu}_{\N^\perp} - \ip{\tau^* \la v_\la}{\tau^* \mu v_\mu}_\L  \\
		&=\ip{w_\la}{w_\mu}_\N + \ip{u_\la}{u_\mu}_{\N^\perp} - ( \ip{w_\la}{w_\mu}_\N + \ip{I(\la)u_\la}{I(\mu)u_\mu}_{\N^\perp}) \\
		&= \ip{(1-I(\mu)^*I(\la))u_\la}{u_\mu}_\M.
\end{align*}
Thus $(\M,u,I)$ is an inner generalized model of $\ph$.

We show next that $\tau$ is a $C$-point for $(\M,u,I)$. To establish this we must produce a vector $u_\tau \in\M$ such that $u_\ln \to u_\tau$ as $n\to\infty$ for every sequence $(\ln)$ in $\D^2$ that converges nontangentially to $\tau$.

As we observed above, $\tau$ is a $B$-point for the model $(\L,v)$ and $\ga\in\ran(1-D\tau)$.    Let $u_\tau$ be the unique element of smallest norm in the nonempty closed convex set $(1-D\tau)^{-1}\ga$.  Then $u_\tau\in\ker(1-D\tau)^\perp=\N^\perp$, and every element of $(1-D\tau)^{-1}\ga$ has the form $e\oplus u_\tau$ for some $e\in\N$.

Let $X_\tau$ be the nontangential cluster set of $v$ at $\tau$ in the model $(\L,v)$; that is, $X_\tau$ comprises the limits in $\L$ of all convergent sequences $(v_\ln)$ for all sequences $(\ln)$ in $\D^2$ that converge nontangentially to $\tau$.  Recall that, by \cite[Proposition 5.8]{AMY10}, a sequence $(v_\ln)$ converges in norm if and only if it converges weakly in $\L$.   If $x\in X_\tau$ is the limit of $u_\ln$ for some sequence $(\ln)$ that converges nontangentially to $\tau$ then, since $(1-D\ln)v_\ln=\ga$, on letting $n\to\infty$ we find that $(1-D\tau)x=\ga$.  Thus
\[
X_\tau \subset (1-D\tau)^{-1}\ga \subset \left\{\begin{pmatrix} e \\ u_\tau \end{pmatrix}: e\in\N\right\}.
\]

We claim that $u_\ln\to u_\tau$ as $n\to\infty$ for every sequence $(\ln)$ in $\D^2$ that converges nontangentially to $\tau$.  For suppose that $u_\ln$ does not converge to $u_\tau$.  Since $v_\ln$, and hence also $u_\ln$, is bounded, on passing to a subsequence we may suppose that $u_\ln \to \xi$ for some vector $\xi\neq u_\tau$, and by passing to a further subsequence, we may suppose that $v_\ln$ converges to some vector $x\in X_\tau$.  But then
\[
v_\ln = \begin{pmatrix} w_\ln \\ u_\ln\end{pmatrix} \to x \in \left\{\begin{pmatrix} e \\ u_\tau \end{pmatrix}: e\in\N\right\},
\]
and hence $u_\ln \to u_\tau$, which is a contradiction.  We have shown that $u_\ln \to u_\tau$ for every sequence $(\ln)$ in $\D^2$ that converges to $\tau$ nontengentially; hence $\tau$ is a $C$-point for the generalized model $(\M,u,I)$.

To see that $\tau$ is a carapoint for $I$, observe that if $\la=r\tau$, where $0<r<1$, then $\la'=(r,r)$, and so by equation \eqref{gotI},
\[
I(r\tau)=1-\frac{(1-r)^2}{1-r}=r.
\]
Hence
\[
\liminf_{\la\to\tau} \frac{1-\|I(\la)\|}{1-\|\la\|_\infty} \leq \liminf_{r\to 1}\frac{1-\|I(r\tau)\|}{1-r} =1.
\]
Thus $\tau$ is a carapoint for $I$.

To complete the proof of condition (3) of Theorem \ref{mainExistence} we must show that $I(\tau) =1_\M$, which by definition means that $I(\la) \to 1_{\N^\perp}$ as $\la \nt \tau$.  Observe that
\beq\label{Iminus1}
I(\la) - 1 =-\frac {(\la'_1-1)(\la'_2-1)}{1-\la'_1(1-Y)-\la'_2Y} =-\bar\tau_1\bar\tau_2\frac{(\la_1-\tau_1)(\la_2-\tau_2)}{1-\bar\tau_1\la_1(1-Y) - \bar\tau_2\la_2 Y}.
\eeq
Since the spectrum of the normal operator $Z=\la'_1(1-Y)+\la'_2Y$ is contained in the convex hull of the points $\la'_1,\la'_2$, 
\begin{align*}
\dist(1,\si(Z)) &\geq \dist(\T,\si(Z)) \geq \dist(\T,\mathrm{conv}\{\la'_1,\la'_2\}) = \dist((\la'_1,\la'_2), \partial \D^2) \\ &= \dist(\la,\partial\D^2).
\end{align*}
It follows that
\[
\norm{(1-\la'_1(1-Y)-\la'_2Y)^{-1}} \leq  \frac{1}{\dist(\la,\partial\D^2)}
\]
and therefore
\[
\norm{I(\la)-1} \leq \frac{|\la_1-\tau_1| \, |\la_2-\tau_2|}{\dist(\la,\partial\D^2)}.
\]
If $\la$ approaches $\tau$ in a set $S$ on which
\[
\frac{\norm{\la-\tau}}{\dist(\la,\partial\D^2)} \leq c < \infty,
\]
then, by the inequality of the means,
\[
\norm{I(\la)-1} \leq \tfrac 12  c \norm{\la-\tau}
\]
for $\la\in S$.  Thus $I(\la) \to 1$ as $\la\nt\tau$.
\end{proof}
A consequence of Theorem \ref{mainExistence} is that $\ph$ has a generalized realization, as in Theorem \ref{genRealiz}.  The preceding proof yields slightly more.
\begin{corollary}
If $\tau\in\T^2$ is a carapoint for $\ph\in\s$ then $\ph$ has a generalized realization
\[
\ph(\la)=a+\ip{I(\la)(1-QI(\la))^{-1}\ga}{\beta}_\M
\]
for some $\beta,\ga \in\M$ and some contraction $Q$ on $\M$ satisfying $\ker (1-Q)=\{0\}$, where $I$ is the inner function given by equation \eqref{defI}, having the properties described in Theorem {\rm \ref{mainExistence}}.
\end{corollary}
\begin{proof}
In the proof of Theorem \ref{mainExistence} it is clear from the definition \eqref{defQ} of $Q$ that $Q$ is a contraction and that $\ker(1-Q)=\{0\}$. From equation \eqref{eq:b2} we have
\[
 \gamma + QI(\lambda)u_\lambda = u_\lambda,
\]
and from the realization $(a, \beta_0, \gamma, D)$ of the model $(\L,v)$, 
 \[
  \vp(\lambda) = a + \ip{\lambda v_\lambda}{\beta_0}.
 \]
Note that, since $\tau^*\beta_0 \in \N^\perp$,
 \begin{align*}
  \vp(\lambda) &= a + \ip{\lambda v_\lambda}{\beta_0}_\L = a+\ip{\tau^*\la v_\la}{ \tau^*\beta_0}_\L       =a + \ip{P_{\N^\perp}\tau^*\la v_\lambda}{\tau^*\beta_0}_{\N^\perp}  \\
 	&= a + \ip{I(\lambda)u_\lambda}{\tau^*\beta_0}_\M.
 \end{align*}
Let $\beta=\tau^*\beta_0 \in \M$.  We then have
\[
 \begin{bmatrix}
  a & 1 \otimes \beta \\
  \gamma \otimes 1 & Q
 \end{bmatrix}
 \vectwo{1}{I(\lambda)u_\lambda} = \vectwo{a + \ip{I(\lambda)u_\lambda}{\beta}}{\gamma + QI(\lambda)u_\lambda} = \vectwo{\vp(\lambda)}{u_\lambda},
\]
and so  $(a, \beta, \gamma, Q)$ is a generalized realization of the generalized model $(\M,u,I)$ of $\ph$.
\end{proof}

We shall call the model $(\M,u,I)$ constructed in the foregoing proof of Theorem \ref{mainExistence} the {\em desingularization} of the model $(\L,v)$ at $\tau$.  The construction depends on the choice of a realization of the model $(\L,v)$, and so where appropriate we should more precisely speak of the desingularization relative to a particular realization.  Of course the singularity of $\ph$ at $\tau$, if there is one, does not disappear; it is shifted into the inner function $I$, where it becomes accessible to analysis by virtue of the formula \eqref{defI} for $I$.

\begin{example}
The inner function $I(.)$ given by equation \eqref{defI} is not in general analytic on $\T^2\setminus \{\tau\}$.

\rm
Let $Y$ be the operation of multiplication by the independent variable $t$ on $L^2(0,1)$ with Lebesgue measure: then $0\leq Y\leq 1$. Let $\tau=(1,1)$.    Suppose that $I$ is analytic at the point $(1,-1)$: then the scalar function 
\[
f(\la) = \ip{I(\la)\mathbf{1}}{\mathbf{1}}
\]
is analytic at $(1,-1)$, where $\mathbf{1}$ denotes the constant function equal to $1$.  We have, for $\la\in\D^2$,
\begin{align*}
f(\la) &= \int_0^1 \frac{t\la_1+(1-t)\la_2-\la_1\la_2}{1-(1-t)\la_1-t\la_2} \ \dd t\\
	&=  \int_0^1 1-\frac{(1-\la_1)(1-\la_2)}{(\la_1-\la_2)t + 1-\la_1} \ \dd t \\
	&= 1- \frac{(1-\la_1)(1-\la_2)}{\la_1-\la_2} [\log((\la_1-\la_2)t+ 1-\la_1)]_0^1 \\
	&= 1- \frac{(1-\la_1)(1-\la_2)}{\la_1-\la_2} [\log(1-\la_2)-\log(1-\la_1)].
\end{align*}
Here we may take any branch of $\log$ that is analytic in $\set{z: \re z > 0}$.
Since $f$ is analytic in a neighborhood of $(1,-1)$, we may let $\la_2 \to -1$ and deduce that,  for some neighborhood $U$ of $1$ and for $\la_1\in U \cap\D$,
\[
f(\la_1, -1) = 1 + 2\frac{1-\la_1}{1+\la_1}[\log 2 - \log(1-\la_1)].
\]
It is then clear that $f(.,-1)$ is not analytic at $1$, contrary to assumption.  Thus $I(.)$ is not analytic at $(1,-1)$, even though $(1,-1) \neq \tau$.
\end{example}

\section{Directional derivatives and slope functions} \label{DirDerivs}
In this section we study the directional derivatives of a function $\ph\in\s$ at a carapoint on the boundary.   One of the main results of \cite{AMY10}, namely Theorem 7.14, asserts the following\footnote{Actually the theorem is slightly more general in that it treats carapoints of $\ph$ in the {\em topological } boundary of $\D^2$.}.
 \begin{theorem}\label{thm2.1}
  Let $\tau \in\T^2$ be a carapoint for $\vp \in \mathcal S$. There exists a function $h$ in the Pick class, analytic and real-valued on $(0, \infty)$, such that the function $z\mapsto -zh(z)$ also belongs to the Pick class,
\beq\label{hof1}
 h(1) = - \liminf_{\lambda \to \tau} \frac{1 - \abs{\vp(\lambda)}}{1 - \norm{\lambda}_\infty} 
\eeq
and, for all $\delta \in \Ha$,
\beq\label{slopeq}
 D_{-\delta}\vp(\tau) = \vp(\tau){\overline{\tau_2}\delta_2}h\left(\frac{\overline{\tau_2}\delta_2}{\overline{\tau_1}\delta_1}\right).
\eeq
 \end{theorem}
With the aid of generalized models we shall present an alternative, more algebraic, proof of this result.  At the same time we obtain further information about directional derivatives at carapoints.
We need a simple  preliminary observation.

\begin{lemma} \label{lem3.1}
  If $\mathcal H$ is a Hilbert space and $Y$ is a positive contraction on $\mathcal H$, then
  \[
  H(z) = -\frac{1}{1-Y+zY}
  \]
  is a well-defined $\mathcal{L(H)}$-valued analytic function on $\C\setminus (-\infty, 0]$. Furthermore, $\IM H(z)$ and $-\IM zH(z)$ are both positive operators for all $z \in \Pi$, and $H(z)$ is Hermitian for $z\in (0,\infty)$.
 \end{lemma}
 \begin{proof}
For any $z\in\C$, the spectrum $\si(1-Y+zY)$ is contained in the convex hull of the points $1, \ z$, by the spectral mapping theorem, and therefore $(1-Y+zY)^{-1}$ is an analytic function of $z$ on the set $\C\setminus (-\infty, 0]$; it clearly takes Hermitian values on the interval $(0,\infty)$. 

   For any $z\in\Pi$ we have
\[
\im (1-Y+zY) = (\im z) Y \geq 0,
\]
and since  $-\im T^{-1}$ is congruent to $\im T$ for any invertible operator $T$, it follows that
\[
\im H(z)= -\im (1-Y+zY)^{-1} \geq 0.
\]
 Similarly
\[
-\im(zH(z))=-\im \frac{1-Y+zY}{z} = -\im \frac{1-Y}{z} = (1-Y)\im\left(- \frac{1}{z}\right) \geq 0.
\]
 \end{proof}

\begin{proof}[Proof of Theorem \rm{\ref{thm2.1}}]
Let $(a,\beta,\ga,D)$ be a realization of $\ph$, associated with a model $(\L,v)$, and let $(M,u,I)$ be the desingularization of this realization at $\tau$.  By Theorem \ref{mainExistence}, $\tau$ is a $C$-point of  $(M,u,I)$, and so there exists $u_\tau\in\M$ such that
\[
\lim_{\la\nt\tau} u_\la = u_\tau   
\]
and, for all $\la,\mu\in\D^2$,
\beq\label{modleq}
 1 - \cc{\vp(\mu)}\vp(\lambda) = \ip{(1 - I\ad(\mu)I(\lambda))u_\lambda}{u_\mu}.
\eeq
Take limits in the last equation as $\mu\nt\tau$ to obtain
\[
 1 - \cc{\vp(\tau)}\vp(\lambda) = \ip{(1 - I(\lambda))u_\lambda}{u_\tau}.
\]
On multiplying through by $-\vp(\tau)$ we deduce that
\begin{align} \label{diffattau}
 \vp(\lambda) - \vp(\tau) &= \vp(\tau)\ip{(I(\lambda) - 1)u_\lambda}{u_\tau} \notag\\
 &= \vp(\tau)\ip{(I(\lambda) - 1)u_\tau}{u_\tau} + \vp(\tau)\ip{(I(\lambda) - 1)(u_\lambda -u_\tau)}{u_\tau}.
\end{align}

Let $\de\in \mathbb{H}(\tau)$, so that $\lambda_t \df \tau - t\delta \in \D^2$ for small enough $t>0$. Then, from equation \eqref{Iminus1},
\begin{align}\label{Ilt-1}
I(\la_t)-1&= I(\tau - t\delta) - 1 = -\bar\tau_1\bar\tau_2\frac{t^2\de_1\de_2}{1-\bar\tau_1(\la_t)_1(1-Y) - \bar\tau_2(\la_t)_2 Y} \nn \\
	&=-\bar\tau_1\bar\tau_2\frac{t\delta_1\delta_2}{\bar\tau_1\delta_1(1 - Y) + \bar\tau_2\delta_2Y}.
\end{align}
In combination with equation \eqref{diffattau} this relation yields
\begin{align*}
 \frac{\vp(\lambda_t) - \vp(\tau)}{t} = -\vp(\tau)&\ip{\frac{\delta_1\delta_2}{\tau_2\delta_1(1 - Y) + \tau_1\delta_2Y}u_\tau}{u_\tau}\\ & - \vp(\tau)\ip{\frac{\delta_1\delta_2}{\tau_2\delta_1(1 - Y) +\tau_1 \delta_2Y}(u_{\lambda_t} -u_\tau)}{u_\tau},
\end{align*}
and on letting $t\to 0+$ we conclude that
 \begin{align*}
D_{-\delta}\vp(\tau) &= -\vp(\tau)\ip{\frac{\delta_1\delta_2}{\tau_2\delta_1(1 - Y) +\tau_1 \delta_2Y}u_\tau}{u_\tau} \\
 &= \vp(\tau)\bar\tau_2\delta_2 h\left(\frac{\bar\tau_2\delta_2}{\bar\tau_1\delta_1}\right)
\end{align*}
where, for any $z\in\Pi$,
\beq\label{defh}
 h(z) = -\ip{\frac{1}{1-Y +z Y}u_\tau}{u_\tau}= \ip{H(z)u_\tau}{u_\tau};
\eeq
here $H(z)$ is as defined in Lemma \ref{lem3.1}.
It is then immediate from Lemma \ref{lem3.1} that $h$ and $-zh(z)$ belong to the Pick class and that $h$  is analytic on $\C\setminus (-\infty, 0]$ and is real-valued on $(0,\infty)$.

It remains to prove equation \eqref{hof1}.  From the definition \eqref{defh} we have
\beq\label{h1w}
h(1)= -\norm{u_\tau}^2,
\eeq
while from the model equation \eqref{modleq}, for any $\la\in\D^2$,
\[
1-|\ph(\la)|^2 = \norm{u_\la}^2 - \norm{I(\la)u_\la}^2.
\]
Let $\la_t= \tau-t\tau$ for $t>0$.  By equation \eqref{Ilt-1} we have
\[
I(\la_t) - 1= -t,
\]
and so, for small enough $t>0$,
\[
1-|\ph(\la_t)|^2 = \norm{u_{\la_t}}^2 - \norm{(1-t)u_{\la_t}}^2 = (2t-t^2)\norm{u_{\la_t}}^2.
\]
We also have
\[
\norm{\la_t}_\infty = \norm{\tau-t\tau}_\infty = (1-t)\norm{\tau}_\infty = 1-t,
\]
and so $1-\norm{\la_t}_\infty^2 = 2t-t^2>0$ for small $t$.  Hence
\[
\frac{1-|\ph(\la_t)|^2}{1-\norm{\la_t}_\infty^2} = \norm{u_{\la_t}}^2
\]
and therefore
\[
\lim_{t\to 0+}\frac{1-|\ph(\la_t)|}{1-\norm{\la_t}_\infty}=\lim_{t\to 0+}\frac{1-|\ph(\la_t)|^2}{1-\norm{\la_t}_\infty^2} = \lim_{t\to 0+}\norm{u_{\la_t}}^2 = \norm{u_\tau}^2.
\]
Hence, by equation \eqref{h1w},
\[
h(1) = -\lim_{t\to 0+}\frac{1-|\ph(\la_t)|}{1-\norm{\la_t}_\infty}.
\]
However, it is known that, for any carapoint $\tau$ of $\ph$,
\[
\lim_{t\to 0+}\frac{1-|\ph(\la_t)|}{1-\norm{\la_t}_\infty} = \liminf_{\la\to\tau} \frac{1-|\ph(\la)|}{1-\norm{\la}_\infty},
\]
(see for example \cite{jaf93} or \cite[Corollary 4.14]{AMY10}).  Equation \eqref{hof1} follows.
\end{proof}

We shall call the function $h$ described in Theorem \ref{thm2.1} the {\em slope function} of $\ph$ at the point $\tau$.  Thus $h$ is the slope function of $\ph$ at a carapoint $\tau\in\T^2$ if, for all $\delta\in\Htau$,
\beq\label{slopeq2}
 D_{-\delta}\vp(\tau) = \vp(\tau){\overline{\tau_2}\delta_2}h\left(\frac{\overline{\tau_2}\delta_2}{\overline{\tau_1}\delta_1}\right).
\eeq
The foregoing proof shows that slope functions have the following representation.
\begin{proposition}
Let $\tau\in\T^2$ be a carapoint for a function $\ph\in\s$.  There exists a Hilbert space $\M$, a vector $u_\tau\in\M$ and a positive contractive operator $Y$ on $\M$ such that, for all $z\in\Pi$,
\beq\label{formSlope}
h(z) = -\ip{\frac{1}{1-Y+zY}u_\tau}{u_\tau}.
\eeq
\end{proposition}

\section{Integral representations of slope functions} \label{represent}
 Theorem \ref{thm2.1} tells us that the directional derivative of a function $\ph\in\s$ at a carapoint is encoded in a slope function $h$ belonging to the Pick class $\Pick$ such that $-zh$ is also in $\Pick$. In this section we derive a representation of functions $h$ with this property.  To obtain such a description we shall need
 the following well-known theorem of Nevanlinna \cite{nev22}, or see \cite[Section II.2, Theorem I]{dono}.

\begin{theorem}\label{thm4.0}
For every holomorphic function $F$ on $\Pi$ such that $\IM F(z) \geq 0$ there exist $c\in\R, \, d\geq 0$ and a finite non-negative Borel measure $\mu$ on $\R$ such that
\beq\label{NevRep}
F(z) = c + dz + \frac{1}{\pi}\int^\infty_{-\infty} \frac{1+tz}{t-z}\dd\mu(t),
\eeq
for all $z\in\Pi$.
Moreover, the $c, d$ and $\mu$ in the representation \eqref{NevRep} are uniquely determined, subject to $c\in\R, \,  d \geq 0, \  \mu\geq 0$ and $\mu(\R) < \infty$.

Conversely, any function $F$ of the form \eqref{NevRep} is in the Pick class.
\end{theorem}
We shall also need another classical theorem -- the Stieltjes Inversion Formula \cite[Section II.2, Lemma I]{dono}.
\begin{theorem}\label{thm4.0.1}
Let $V$ be a nonnegative harmonic function on $\Pi$, and suppose that $V$ is the Poisson integral of a positive measure $\mu$ on $\R$:
\beq\label{poisson}
V(x+\ii y) = cy+ \frac{y}{\pi}\int^\infty_{-\infty} \frac{\dd\mu(t)}{(t-x)^2+y^2}
\eeq
for some $c\geq 0$ and all $y>0$, where 
\beq\label{muflfin}
\int^\infty_{-\infty} \frac{\dd\mu(t)}{1+t^2} < \infty.
\eeq
Then
\begin{equation}
\lim_{y \to 0+} \int_a^b V(x + \ii y)\ \dd x = \mu((a,b)) + \tfrac 12 \mu(\{a\}) + \tfrac 12 \mu(\{b\})
\end{equation}
whenever $-\infty < a < b < \infty$.
\end{theorem}
We can now identify the class of $h \in \mathcal P$ such that  $-zh \in \mathcal P$.
 \begin{theorem}\label{thm4.1}
  The following are equivalent for any analytic function $h$ on $\Pi$.
\begin{enumerate}
 \item[\rm (i)] $h, -zh \in \mathcal P$;
 \item[(ii)] $h \in \mathcal P$ and the Nevanlinna representation of $h$ has the form
\[
 h(z) = c + dz + \frac{1}{\pi}\int\frac{1 + tz}{t - z}\dd\mu(t)
\]
where
\begin{enumerate}
\item $d = 0$,
 \item $\mu((0,\infty)) = 0$  and
 \item $\ds c \leq \frac{1}{\pi}\int t \ \dd\mu(t)$;
\end{enumerate}
 \item[(iii)] there exists a positive Borel measure $\nu$ on $[0,1]$ such that 
\[
 h(z) = -\int\frac{1}{1 - s + sz} \dd\nu(s).
\]
\end{enumerate}
 \end{theorem}

\begin{proof}
 (i)$\Rightarrow$(ii) Let $h$ and $-zh$ be in the Pick class.  Then there exist unique $c, c' \in \R, d, d' \geq 0$, and finite positive Borel measures $\mu, \ \nu$ on $\R$ such that
\beq\label{hz}
 h(z) = c + dz + \frac{1}{\pi}\int\frac{1 + tz}{t-z}\dd\mu(t)
\eeq
and
\beq\label{zhz}
 -zh(z) = c' + d'z + \frac{1}{\pi}\int\frac{1 + tz}{t-z}\dd\nu(t).
\eeq
If $z=x+\ii y$, where $x,y\in\R$, then
\[ 
\IM\frac{1 + tz}{t-z} = \im \left(\frac{1+t^2}{t-z} -t\right)=(1+t^2)\IM\frac{1}{t-z}
 = (1+t^2)\frac{y}{(t-x)^2 + y^2}.
\]
Hence
\begin{align*}
 \IM h(z) = dy + \frac{y}{\pi}\int\frac{1}{(t-x)^2 + y^2}(1+t^2)\dd\mu(t), \\
 \IM (-zh(z)) = d'y + \frac{y}{\pi}\int\frac{1}{(t-x)^2 + y^2}(1+t^2)\dd\nu(t).
\end{align*}
Since $\IM h$ is nonnegative and harmonic, Theorem \ref{thm4.0.1} implies that
\begin{equation}
 \lim_{y\to 0^+} \int^b_a \IM h(x+iy)\dd x = \mu((a,b)) + \frac{\mu(\{a\}) + \mu(\{b\})}{2} \label{eq4.0}
\end{equation}
and
\begin{equation}
 \lim_{y \to 0^+} \int^b_a \IM(-zh(z))\dd x = \nu((a,b)) + \frac{\nu(\{a\}) + \nu(\{b\})}{2}. \label{eq4.0.1}
\end{equation}
Note that 
\begin{equation}
\IM(-zh) = -\IM((x+iy)h) = -x\IM h - y \RE h, \label{eq4.1}
\end{equation} 
and so
\begin{equation}
 \lim_{y \to 0^+} \int^b_a \IM(-zh(z))\ \dd x = -\lim_{y \to 0^+}\int^b_a x\IM h(x+iy)\ \dd x - \lim_{y \to 0^+} y\int^b_a \RE h(x+iy)\ \dd x \label{eq4.2}.
 \end{equation}
 Now let 
 \[
 A_y = \int^b_a x\IM h(x+iy)\ \dd x 
 \]
 and
 \[
 B_y = y\int^b_a \RE h(x+iy)\ \dd x,
 \]
 so that 
 \begin{equation}\label{AyBy}
 \lim_{y \to 0^+} \int^b_a \IM(-zh(z))\ \dd x =  -\lim_{y \to 0^+}A_y  - \lim_{y \to 0^+} B_y.
\end{equation}

\begin{lemma}\label{limBy}
For any $a,b\in\R$ such that $a<b$,
\[
 \lim_{y \to 0^+} B_y = \lim_{y \to 0^+} y\int^b_a \RE h(z)\ \dd x = 0.
\]
\end{lemma}
\begin{proof}
In view of the representation \eqref{hz} of $h$ we have
\begin{align}\label{By}
 B_y &= y\int_a^b \re \int\frac{1+t(x+\ii y)}{t-x-\ii y} \ \dd\mu(t) \  \dd x   \nn \\
	&= y\int_a^b \re \int\frac{(1+t(x+\ii y))(t-x+\ii y)}{(t-x)^2 + y^2} \ \dd\mu(t) \  \dd x  \nn  \\
	&=y\int_a^b  \int\frac{(1+tx)(t-x)- t y^2}{(t-x)^2 + y^2} \ \dd\mu(t) \  \dd x  \nn \\
	&=y\int_a^b  \int\frac{(t-x)(1+(t-x)x+x^2)- (t-x) y^2-xy^2}{(t-x)^2 + y^2} \ \dd\mu(t) \  \dd x   \nn\\
	&=y\int_a^b  \int xC_2 + (1+x^2-y^2) C_1 - xy^2 C_0 \ \dd\mu(t) \  \dd x  
\end{align}
where
\[
C_2= \frac{(t-x)^2}{(t-x)^2+y^2}, \quad C_1=\frac{t-x}{(t-x)^2+y^2} \quad\mbox{ and } \quad C_0=\frac{1}{(t-x)^2+y^2}.
\]
For all $t, x$ in $\R$ and $y>0$  we have $C_2 \leq 1$ and $C_0 \leq 1/y^2$, and so
\beq\label{estC02}
\left|  xC_2  - xy^2 C_0\right| \leq  2|x|.
\eeq
Choose $N \geq 1+\max\{\abs{a},\abs{b}\}$ such that $\mu(\{N,-N\})=0$.  Then
\begin{align}\label{estC20}
\int_a^b  \int |xC_2  - xy^2 C_0| \ \dd\mu(t) \  \dd x &\leq \int_a^b \int 2|x| \ \dd\mu(t) \  \dd x \nn \\
	&\leq\mu(\R)\int_{-N}^N 2|x|   \dd x \nn \\
	&=2N^2\mu(\R).
\end{align}
It is then immediate that
\beq\label{limC20}
\lim_{y\to 0+} y\int_a^b  \int |xC_2  - xy^2 C_0| \ \dd\mu(t) \  \dd x = 0.
\eeq

For $|t|\geq N, \ a\leq x\leq b$ we have $|t-x| \geq 1$,  hence $|C_1|\leq 1$ and so
\begin{align}\label{C1gtN}
\int_a^b  \int_{|t|\geq N} |(1+x^2+y^2) C_1| \ \dd\mu(t) \  \dd x &\leq \int_a^b  \int 1+x^2+y^2 \ \dd\mu(t) \  \dd x  \nn \\
	&\leq \mu(\R)(1+N^2+y^2)(b-a).
\end{align}
On the other hand, when $|t| \leq N$ and $a\leq x \leq b$,
\[
| (1+x^2+y^2) C_1| \leq(1+N^2+y^2) \frac{ |t-x|}{ (t-x)^2 + y^2}.
\]
On making the change of variable $s=|t-x|$  and observing that $ 0 \leq s \leq 2N$ when $|t| \leq N$ and $a \leq x \leq b$, we find that
\begin{align*}
\int_a^b   |(1+x^2+y^2) C_1 |  \  \dd x  &\leq  2(1+N^2+y^2) \int_0^{2N} \frac{s \ \dd s}{s^2+y^2} \\
	&= (1+N^2+y^2) \left( \log(4N^2+y^2)-2\log y\right),
\end{align*}
and therefore
\beq\label{C1ltN}
\int_{|t|\leq N} \ \dd\mu(t) \int_a^b   |(1+x^2+y^2) C_1 |  \  \dd x \leq \mu(\R) (1+N^2+y^2) \left( \log(4N^2+y^2)-2\log y\right) < \infty.
\eeq
It follows from the Fubini-Tonelli theorem that the order of integration can be reversed, and on combining the estimates \eqref{C1gtN} and \eqref{C1ltN} we find that 
\[
\int_a^b  \int |(1+x^2+y^2) C_1| \ \dd\mu(t) \  \dd x \leq \mu(\R)(1+N^2+y^2)\left[b-a+  \log(4N^2+y^2)-2\log y\right],
\]
from which it is clear that
\[
\lim_{y\to 0+} y\int_a^b  \int |(1+x^2+y^2) C_1| \ \dd\mu(t) \  \dd x =0.
\]
On combining this statement with \eqref{limC20} we conclude that 
\[
\lim_{y\to 0+} y\int_a^b  \int |xC_2+(1+x^2+y^2) C_1-xy^2C_0| \ \dd\mu(t) \  \dd x =0
\]
and hence, by equation \eqref{By},  that $B_y \to 0$ as $y\to 0+$.
\end{proof}

Since $\im h \geq 0$ we have
\[
 a\int^b_a \IM h(x+iy)\ \dd x \leq A_y \leq  b\int^b_a \IM h(x+iy) \dd x.
\]
Combining this inequality with (\ref{eq4.0}) we find that
\[
 a\left(\mu(a,b) +\frac{\mu(\{a\}) + \mu(\{b\})}{2}\right) \leq \lim_{y \to 0^+} A_y \leq  b\left(\mu(a,b) +\frac{\mu(\{a\}) + \mu(\{b\})}{2}\right)
\]
and so, in view of equation (\ref{eq4.0.1}), for all $a < b$,
\begin{align} \label{mu&nu}
-b\left(\mu(a,b) +\frac{\mu(\{a\}) + \mu(\{b\})}{2}\right) &\leq \nu((a,b)) + \frac{\nu(\{a\}) + \nu(\{b\})}{2} \\
  & \leq  -a\left(\mu(a,b) +\frac{\mu(\{a\}) + \mu(\{b\})}{2}\right).
\end{align}
As this inequality holds for all $a < b \in \R$, we can let $a = 0$ and $b >0$. Then 
\[
 \nu((0,b)) + \frac{\nu(\{0\}) + \nu(\{b\})}{2} \leq  0.
\]
But as $\nu$ is a positive measure, this implies that $\nu((0, \infty)) = 0$ and $\nu(\{0\}) = 0$, i.e. $\nu([0,\infty)) = 0$.

Now let $0<a<b$. Then
\[
 a\left(\mu(a,b) +\frac{\mu(\{a\}) + \mu(\{b\})}{2}\right) \leq \left(\nu(a,b) +\frac{\nu(\{a\}) + \nu(\{b\})}{2}\right) = 0.
\]
But $\mu \geq 0$, and so $\mu((a,b)) = 0$.  It follows that $\mu((0,\infty)) = 0$, which is to say that condition (b) holds. 

\begin{fact*}
 For $t<0, \nu(\{t\}) = -t \mu(\{t\})$. 
\end{fact*}
 Since $\mu, \nu$ are finite and positive, they can only have at most countably many point masses, and so we may choose a sequence of intervals $(a_n,b_n)\subset (-\infty, 0)$ such that $\mu(\{a_n\}) = \mu(\{b_n\}) = \nu(\{a_n\}) = \nu(\{b_n\}) = 0$,  $t \in (a_n, b_n)$ for all $n$  and $\bigcap (a_n, b_n) = \{t\}$.   Inequality \eqref{mu&nu} implies that
\[
 -b_n\mu((a_n, b_n)) \leq \nu((a_n, b_n)) \leq -a_n\mu((a_n, b_n)),
\]
and in the limit we obtain
\[
 -t\mu(\{t\}) = \nu(\{t\}) \leq -t\mu(\{t\}).
\]

If $\sigma$ is a finite positive measure on $(-\infty,0)$, we shall  call a finite partition $P= \{x_1,\dots,x_n\}$, where $ x_1 < x_2 < ... < x_n < 0$,  \emph{special for $\si$} if $\sigma(P) =0$. 

\begin{fact*}
 If $f$ is continuous on $(-\infty,0)$ with compact support and $\epsilon > 0$, there exists a partition $P$ that is special for $\sigma$ such that
 \[
 \abs{\int f  \ \dd\sigma - S(f, P)} < \epsilon,
 \]
 where $S(f,P)$ denotes the Riemann sum of $f$ over $P$.
\end{fact*}

\begin{lemma}\label{lem5.5}
 If $f$ is a continuous function of compact support on $(-\infty, 0)$ then 
 \[
\int f \  \dd\mu  = \int f (t) \  t\dd\nu(t).
 \] 
\end{lemma}

Equations \eqref{hz} and \eqref{zhz} give us two different expressions for $-zh(z)$:
\[
 -z\left(c +dz +\frac{1}{\pi}\int \frac{1 +tz} {t -z}\ \dd\mu(t)\right) = c' + d'z + \frac{1}{\pi} \int \frac{1 + tz}{t - z}\ \dd\nu(t).
 \]
 Hence, by Lemma \ref{lem5.5},
 \[
 -z\left(c +dz +\frac{1}{\pi}\int \frac{1 +tz} {t -z} \ \dd\mu(t)\right) = c' + d'z + \frac{1}{\pi} \int \frac{1+tz}{t - z} \ (-t\dd\mu(t)).
\]
We may rearrange this equation, in order to compare polynomials, obtaining
\begin{align*}
 c' + (d' + c)z +dz^2 &= \frac{1}{\pi} \int \frac{1 + tz}{t - z}(t-z) \ \dd\mu(t) \\
 &= \frac{1}{\pi} \int (1+tz) \ \dd\mu(t) \\
 &= \frac{1}{\pi}\int \ \dd\mu(t) + \left(\frac{1}{\pi}\int t\dd\mu(t)\right)z.
\end{align*}
We immediately see that 
\[
c' = \frac{1}{\pi}\int \ \dd\mu, \quad  d'+c = \frac{1}{\pi}\int  t \ \dd\mu(t) \quad \mbox{ and } \quad  d = 0.
\]
The last of these statements is condition (a) in (ii).  Since $d' > 0$ the second statement tells us that
\[
 c \leq d' + c = \int t \ \dd\mu(t),
\]
which is condition (c).  This concludes the proof that (i)$\Rightarrow$(ii). 

(ii)$\Rightarrow$(iii) 
Suppose that $h\in\Pick$ has a Nevanlinna representation that satisfies conditions (a)-(c) of (ii).  Then, for $z\in\Pi$,
\begin{align}\label{eq4.3}
 h(z) &= c + \frac{1}{\pi}\int\frac{1+tz}{t-z} \ \dd\mu(t) \notag \\
 &= c + \frac{1}{\pi}\int \left( \frac{1+t^2}{t - z}-t\right) \ \dd\mu(t) \notag \\
 &= c - \frac{1}{\pi}\int t\ \dd\mu(t) + \frac{1}{\pi}\int\frac{1 - t}{t - z}\frac{1+t^2}{1 - t}\ \dd\mu(t).
\end{align}

Since the indefinite integral of $\frac{1 + t^2}{1 - t}\ \dd\mu(t)$ is a finite positive measure on $(-\infty, 0]$, we may define a finite positive Borel measure $\nu$ on $[0,1]$ by 
\begin{align}
\nu(\{0\}) &= \frac{1}{\pi}\int t\ \dd\mu(t) - c, \label{eq4.4} \\
\nu(E) &= \frac{1}{\pi}\int_{\tilde E} \frac{1 + t^2}{1 - t}\ \dd\mu(t) \label{eq4.5}
\end{align}
for any Borel set $E\subset (0,1]$, where $\tilde E \df \{1-1/s: s\in E\}$.

With this definition, if $\psi$ is a continuous bounded function on $(-\infty, 0]$,
\[
 \frac{1}{\pi}\int \psi(t) \frac{1 + t^2}{1 - t} \ \dd\mu(t) = \int_{(0,1]} \psi\left(1 - \frac{1}{s}\right)\ \dd\nu(s).
\]

From equations (\ref{eq4.3}) and (\ref{eq4.5}),
\begin{align*}
 h(z) &= c - \frac{1}{\pi}\int t \ \dd\mu(t) + \frac{1}{\pi}\int \left(\frac{1-t}{t-z}\right)\frac{1 + t^2}{1 - t}\ \dd\mu(t) \\
 &= -\nu(\{0\}) + \int_{(0,1]}\frac{1 - (1-\frac{1}{s})}{1 - \frac{1}{s} -z} \ \dd\nu(s) \\
 &= -\nu(\{0\}) + \int_{(0,1]}\frac{1}{s - 1 - sz} \ \dd\nu(s) \\
 &= -\left[\nu(\{0\}) + \int_{(0,1]}\frac{1}{1 - s + sz} \ \dd\nu(s)\right] \\
 &= -\int_{[0,1]}\frac{1}{1 - s + sz} \ \dd\nu(s), 
\end{align*}
which completes the proof that (ii)$\Rightarrow$(iii).

(iii)$\Rightarrow$(i) Suppose that $\nu$ is a positive finite Borel measure  on $[0,1]$ and 
\[
h(z) = -\int \frac{1}{1 - s + sz}\ \dd\nu(s)
\]
for all $z\in\Pi$.
Let $Y$ be the operator of multiplication by the independent variable $s$ on $L^2(\nu)$. Evidently $Y$ is a positive contraction, and hence,  by Lemma \ref{lem3.1},  for any $z\in\Pi$, the operators
\[
  -\im (1-Y+Yz)^{-1} \quad \mbox{  and  } \quad  \im \left( z(1 - Y + Yz)^{-1} \right)
\]
on $L^2(\nu)$ are positive definite.  Since 
\[
 \im h(z) = -\im \int\frac{1}{1-s+sz}\ \dd\nu(s) = \ip{-\im\frac{1}{1-Y+Yz}1}{1}_{L^2(\nu)} \geq 0
\]
and likewise
\[
\im (-zh(z)) = \im \int\frac{z}{1-s+sz}\ \dd\nu(s) = \ip{\im\frac{z}{1-Y+Yz}1}{1}_{L^2(\nu)} \geq 0,
\]
it follows that (i) holds. 
\end{proof}
The proof shows that if $h$ and $-zh$ belong to $\Pick$ then $h$ is analytic on $(0,\infty)$.

\section{Functions with prescribed slope function} \label{converse}

   In this section we prove a converse to  Theorem \ref{thm2.1}:  we construct, for any function $h\in \Pick$ such that $-zh \in\Pick$, a function $\ph\in\s$ with slope function $h$ at a carapoint.

We shall need the following simple observation about the Cayley transform (an aplication of the quotient rule).   The {\em two-variable Herglotz class}  is defined to be the set of analytic functions on $\D^2$ with non-negative real part.
\begin{lemma}\label{lem5.1}
If $f$ is a  function in the two-variable Herglotz class then the function $\vp$ on $\D^2$ given by
\[
\vp = \frac{1-f}{1+f}
\]
belongs to $\s$. Furthermore,  if $\tau\in \T^2$ is such that the radial limit 
\[
f(\tau) \df \lim_{r\to 1-} f(r\tau)
\]
exists and is not $ -1$ and if the directional derivative $D_{-\delta} f(\tau)$ exists for some direction $\de$, then so does $D_{-\delta} \vp (\tau)$, and
\begin{equation} \label{eq5.1}
  D_{-\delta} \vp (\tau) = \frac{-2 D_{-\delta} f(\tau)}{(1+f(\tau))^2}. 
\end{equation}
\end{lemma}
Recall that $\chi$ denotes the point $(1,1)$.
\begin{theorem}\label{SuffSlope}
If $h \in \mathcal P$  and $-zh \in \mathcal P$ then there exists $\vp \in \s$ such that $\chi$ is a carapoint for $\vp$, $\ph(\chi)=1$ and $h$ is the slope function for $\ph$ at $\chi$.
\end{theorem}
\begin{proof}
   By Theorem \ref{thm2.1} there exists a positive Borel measure $\nu$ on $[0,1]$ such that 
 \[
  h(z) = -\int\frac{1}{1-s +sz}\ \dd\nu(s).
 \]
 Define a family of functions $f_s$ on $\D^2$ for $s \in [0,1]$ by
 \[
  f_s(\la)  = \left(s\frac{1+\lambda_1}{1-\lambda_1} + (1-s)\frac{1+\lambda_2}{1-\lambda_2}\right)^{-1}.
 \]
  For any $\la\in\D^2$ the denominator on the right hand side is a convex combination of two points in $\Ha$, hence belongs to $\Ha$.  Thus $f_s$ lies in the two-variable Herglotz class for $0\leq s\leq 1$.  Moreover, for $0<r<1$ and every $s\in[0,1]$,
\beq\label{valfrchi}
f_s(r\chi) = \frac{1-r}{1+r}
\eeq
and hence the radial limit
\[
f_s(\chi) \df \lim_{r\to 1-} f_s(r\chi) = \lim_{r\to 1-} \frac{1-r}{1+r}
\]
exists and is zero.
 We compute $D_{-\delta} f_s(\chi)$. 
 \begin{align}\label{diffQuo}
  \frac{f_s(\chi - t\delta) - f(\chi)}{t} &=\frac{1}{t} \left(s\frac{1 + 1 - t\delta_1}{1 - (1-t\delta_1)} + (1-s)\frac{1 + 1 - t\delta_2}{1 - (1-t\delta_2)}\right)^{-1}  \nn \\
  &= \left(s\frac{2}{\delta_1}+(1-s)\frac{2}{\delta_2}-t\right)^{-1} \nn \\
	&= \frac{1}{2} \frac{\de_1\de_2}{(1-s)\de_1 + s\de_2-\half t \de_1\de_2}. 
 \end{align}
On letting $t \to 0$ we obtain
 \beq\label{Ddefs}
  D_{-\delta} f_s(\chi) = \frac{1}{2} \frac{\de_1\de_2}{(1-s)\delta_1+ s\delta_2}.
 \eeq
 
Define a function $f$ on $\D^2$ by
 \[
  f(\la) = \int f_s(\la) \  \dd\nu(s).
 \]
Since $\re f_s(\la) > 0$ for every $s\in  [0,1]$, it is clear that $f$ lies in the two-variable Herglotz class.  Furthermore, in view of equation \eqref{valfrchi}, for $0<r<1$,
\beq\label{frchi}
f(r\chi) = \frac{1-r}{1+r} \nu[0,1]
\eeq
and $f$ has radial limit $0$ at $\chi$: 
\beq\label{radial}
f(\chi) \df \lim_{r\to 1-} \int f_s(r\chi) \ \dd\nu(s) =  \lim_{r\to 1-}\int  \frac{1-r}{1+r} \ \dd\nu(s) =  \lim_{r\to 1-} \nu[0,1] \frac{1-r}{1+r} = 0.
\eeq
  Let us calculate the directional derivative of $f$ at $\chi$ in the direction $-\de$ where $\de\in\HH$.  Equation \eqref{Ddefs} suggests that
  \begin{equation}\label{eq5.3}
  D_{-\delta}f(\chi) = \frac{1}{2}\int\frac{\de_1\delta_2}{(1-s)\de_1 + s\delta_2}\ \dd\nu(s). 
 \end{equation}
We must verify that this is correct.  By equation \eqref{diffQuo}, we have, for small $t>0$,
\begin{align}\label{estRem}
\frac{f(\chi-t\de)-f(\chi)}{t} &-\frac{1}{2}\int\frac{\de_1\delta_2}{(1-s)\de_1 + s\delta_2}\ \dd\nu(s)\nn \\
	&= \frac{1}{2} \int\frac{\de_1\de_2}{(1-s)\de_1 + s\de_2-\half t \de_1\de_2} - \frac{\de_1\de_2}{(1-s)\de_1 + s\de_2} \ \dd\nu(s)\nn  \\
	&= \frac{\de_1\de_2}{2} \int \frac{\half t\de_1\de_2 \ \dd\nu(s)}{\left((1-s)\de_1+s\de_2-\half t\de_1\de_2\right)\left((1-s)\de_1+s\de_2\right)}.
\end{align}
Since $\de_1,\de_2\in\Ha$, the distance $K$ from $0$ to the convex hull of $\{\de_1,\de_2\}$ is positive.  For sufficiently small $t>0$ we have, for all $s\in[0,1]$,
\[
|(1-s)\de_1+s\de_2-\half t\de_1\de_2| \geq \half K,
\]
and for such $t$ the denominator of the integrand in equation \eqref{estRem} is at least $\half K^2$.  It follows that, for small enough $t$,
\[
\left|\frac{f(\chi-t\de)-f(\chi)}{t} -\frac{1}{2}\int\frac{\de_1\delta_2}{(1-s)\de_1 + s\delta_2}\ \dd\nu(s)\right| \leq\frac{|\de_1\de_2|^2\nu[0,1]}{2K^2}t,
\]
and hence equation \eqref{eq5.3} is correct.

Let $\vp$ be defined by
 \[
  \vp = \frac{1 - f}{1 + f}.
 \]
We claim that $\chi$ is a carapoint for $\ph$.  For any $\la\in\D^2$,
\[
\frac{1-|\ph(\la)|^2}{1-\|\la\|_\infty^2} = \frac{4\re f(\la)}{(1-\|\la\|_\infty^2)|1+f(\la)|^2}
\]
and so, by equation \eqref{frchi},
\begin{align*}
\frac{1-|\ph(r\chi)|^2}{1-\|r\chi\|_\infty^2} &= \frac{4\nu[0,1]}{\left(1+r+(1-r)\nu[0,1]\right)^2} \\
	& \to \nu[0,1] \quad \mbox{ as } r \to 1-.
\end{align*}
  Hence
\[
\liminf_{\la\to \chi} \frac{1-|\ph(\la)|^2}{1-\|\la\|_\infty^2} \leq \nu[0,1] < \infty
\]
and $\chi$ is a carapoint for $\ph$.

 In view of equation \eqref{radial} it is clear that $\ph$ has radial limit $1$ at $\chi$, that is to say, $\ph(\chi)=1$.  By Lemma \ref{lem5.1},  $\vp$ lies in $\s$ and has directional derivative at $\chi$ given by
 \begin{align*}
  D_{-\delta} \vp(\chi) &= \frac{-2D_{-\delta} f(\chi)}{(1+f(\chi))^2} \\
  &= (-2) \frac{1}{2} \int \frac{\de_1\de_2}{(1-s)\de_1 + s\delta_2}\ \dd\nu(s) \\
  &= -\delta_2h\left(\frac{\delta_2}{\delta_1}\right)\\
	&= -\ph(\chi) \delta_2h\left(\frac{\delta_2}{\delta_1}\right).
 \end{align*}
Thus $h$ is the slope function for $\ph\in\s$ at the point $\chi$.
\end{proof}
By a simple change of variable we obtain the following.
\begin{corollary}\label{prescription}
Let $\omega\in\T$, let $\tau\in\T^2$ and let $h,\ -zh(z) \in\Pick$.  There exists a function $\ph\in\s$ having a carapoint at $\tau$ such that $\ph(\tau)=\omega$ and $h$ is the slope function of $\ph$ at $\tau$.
\end{corollary}

\section{Nevanlinna representations in two variables}\label{NevanRep}

The following refinement of Theorem \ref{thm4.0}, also due to Nevanlinna, is the main tool in one of the standard  proofs of the Spectral Theorem for unbounded self-adjoint operators \cite{lax}.
\begin{proposition}\label{Nev2}
Let $h \in \mathcal P$. If 
\beq\label{NevCond}
\lim_{y\to\infty}y \IM h(iy) < \infty
\eeq
 then there exists a finite  positive measure $\mu$ on $\R$ such that, for all $z\in\Pi$,
\beq\label{NevRep2}
 h(z) = \int\frac{\dd\mu(t)}{t-z}\ .
\eeq
\end{proposition}
For a proof see \cite{lax}.   

In this section we shall generalize Proposition \ref{Nev2} to two variables.  We need an  analog for the Cauchy transform formula \eqref{NevRep2}.  The closest one we can find involves the two-variable resolvent of a self-adjoint operator $B$ on a Hilbert space $\M$, to wit
\[
h(z_1,z_2) = b-\ip{\left(B+z_1 Y + z_2(1-Y)\right)^{-1}\al}{\al}
\]
for some $b\in\R, \ \al\in\M$ and some positive contraction $Y$ on $\M$.  In an earlier paper \cite[Theorem 6.9]{AMY11} we obtained a somewhat similar result, but with the unsatisfactory feature that the representation obtained was not of $h$ itself but rather of a ``twist" of $h$.  The use of generalized models enables us to remedy this defect.

The growth condition \eqref{NevCond} is expressible in terms of carapoints of the Schur-class function $\ph$ associated with $h$ by the definition
\beq\label{ph&h}
\ph(\la)= \frac{h(z)-i}{h(z)+i}\quad  \mbox{ where } \quad z= i\frac{1+\la}{1-\la}.
\eeq
Let us establish the corresponding assertion for functions of two variables.  We denote by $\p$ the two-variable Pick class, that is the set of analytic functions on $\Pi^2$ with non-negative imaginary part and we recall that $\chi$ denotes $(1,1)$.
\begin{proposition}\label{carainfty}
Let $h\in\p$ and let $\ph\in\s$ be defined by
\beq\label{phiandh}
\ph(\la) = \frac{h(z)-i}{h(z)+i} \quad \mbox{ where } z_j= i\frac{1+\la_j}{1-\la_j}, \quad j=1,2,
\eeq
for $\la\in\Pi^2$.  The following conditions are equivalent.
\begin{enumerate}
\item $\liminf_{y\to \infty} y\im h(iy\chi) < \infty$;
\item $\lim_{y\to\infty} y\im h(iy\chi)$ exists and is finite;
\item $\chi$ is a carapoint for $\ph$ and $\ph(\chi)\neq 1$;
\item $(0,0)$ is a carapoint for the function $H\in\p$ given by $H(z)= h(-1/z_1, -1/z_2)$.
\end{enumerate}
\end{proposition}
\begin{proof}
(2)$\Rightarrow$(1) is trivial.

(1)$\Rightarrow$(3)  Suppose (1) holds and let $\beta$ be the limit inferior in (1).  There is a sequence $(y_n)$ in $\R^+$ such that $y_n \to \infty$ and
\[
\lim_{n\to \infty} y_n \im h(iy_n\chi) = \beta.
\]
Let
\[
r_n= \frac {iy_n-i}{iy_n+i} = \frac{y_n-1}{y_n+1}.
\]
Then $y_n = \frac{1+r_n}{1-r_n}$ and $r_n \to 1-$ as $n\to\infty$.  From the relation
\[
1-|\ph(\la)|^2 = \frac{4\im h(z)}{|h(z)+i|^2}
\]
we have
\[
1-|\ph(r_n\chi)|^2 = \frac{4\im h(iy_n\chi)}{|h(iy_n\chi)+i|^2}.
\]
Since $|h(z)+i| \geq 1$ for all $z\in\Pi^2$,
\[
1-|\ph(r_n\chi)|^2 \leq 4\im h(iy_n\chi).
\]
Similarly
\begin{align*}
1-\|r_n\chi\|_\infty^2 &= 1-r_n^2 = \frac{4\im iy_n}{|iy_n+i|^2 }\\
	&=\frac{4y_n}{(1+y_n)^2}.
\end{align*}
Hence
\begin{align*}
\frac{1-|\ph(r_n\chi)|^2}{1-\|r_n\chi\|_\infty^2} &\leq 4\im h(iy_n\chi)\frac{(1+y_n)^2}{4y_n} \\
	&\to \beta \quad \mbox{ as } n\to \infty.
\end{align*}
Consequently
\[
\liminf_{\la\to\chi} \frac{1-|\ph(\la)|^2}{1-\|\la\|_\infty^2} \leq \beta <\infty,
\]
and so $\chi$ is a carapoint for $\ph$.

By the Carath\'eodory-Julia theorem for the bidisc \cite{jaf93,AMY10},
\beq\label{defalpha}
 \al \df \liminf_{\la\to\chi} \frac{1-|\ph(\la)|^2}{1-\|\la\|_\infty^2} = \lim_{r\to 1-}\frac{1-|\ph(r\chi)|^2}{1-r^2}> 0
\eeq
and $\beta\neq 0$ since $\al\leq\beta$.
Now for any $r\in(0,1), \, y=\frac{1+r}{1-r}$, a simple calculation shows that
\begin{align}\label{yimhychi}
y\im h(iy\chi) &= \frac{1+r}{1-r} \, \frac{1-|\ph(r\chi)|^2}{|1-\ph(r\chi)|^2} \nn \\
	&=\frac{(1+r)^2}{|1-\ph(r\chi)|^2} \, \frac{1-|\ph(r\chi)|^2}{1-r^2}.
\end{align}
On putting $r=r_n$ and letting $n\to\infty$ we find that 
\[
\lim_{n\to\infty} |1-\ph(r_n\chi)|^2 = \frac{4\al}{\beta} \neq 0.
\]
Thus $\ph(\chi) \neq 1$.   Hence (1)$\Rightarrow$(3).

(3)$\Rightarrow$(2)  Suppose (3).  Then the quantity $\al$ defined by equation \eqref{defalpha} satisfies $0< \al <\infty$.  On letting $y\to \infty$ (and hence $r\to 1-$) in equation \eqref{yimhychi} we obtain
\[
\lim_{y\to\infty} y \im h(iy\chi) = \frac{4\al}{|1-\ph(\chi)|^2},
\]
and so (2) holds.

(2)$\Leftrightarrow$(4) According to the general definition of a carapoint in Section \ref{intro}, $(0,0)$ is a carapoint for $H\in\p$ if 
\[
\liminf_{z\to(0,0)} \frac{\im H(z)}{\min\{\im z_1,\im z_2\}} < \infty,
\]
and by the two-variable Carath\'eodory-Julia theorem \cite{jaf93,AMY10}, this is so if and only if
\[
\liminf_{\eta\to 0} \frac{\im H(i\eta\chi)}{\eta}  =\liminf_{\eta\to 0} \frac{\im h(i\chi/\eta)}{\eta}< \infty.
\]
On setting $y=1/\eta$ we deduce that (2)$\Leftrightarrow$(4).
\end{proof}
We shall say that $\infty$ is a {\em carapoint for $h\in\p$ with finite value} if the equivalent conditions of Proposition \ref{carainfty} hold.  We define the value $h(\infty)$ in this case by
\[
h(\infty) = \lim_{y\to \infty} h(iy\chi) = H(0,0) = i\frac{1+\ph(\chi)}{1-\ph(\chi)}
\]
where $H, \, \ph$ are as in Proposition \ref{carainfty}.  There is also a notion of carapoint of $h$ with infinite value: see \cite[Section 7]{ATY}.

Here is our generalization of the Nevanlinna representation \eqref{NevRep2} to the two-variable Pick class.
\begin{theorem}\label{NevanlinnaRep}
 The following statements are equivalent for a function $h:\Pi^2\to\C$.   
\begin{enumerate}
\item  $h$ is in the Pick class $\p$, $\infty$ is a carapoint for $h$ with finite value;
\item there exist a scalar $b\in \R$, a Hilbert space $\M$, a vector $\al\in\M$, a positive contraction $Y$ on $\M$ and a densely defined self-adjoint operator $B$ on $\M$ such that, for all $z\in\Pi^2$,
  \beq\label{RepForm}
   h(z) = b - \ip{\left(B + z_1Y + z_2(1-Y)\right)^{-1}\alpha}{\alpha}.
  \eeq
\end{enumerate}
\end{theorem}
\begin{proof}
We begin by observing that the inverse in equation \eqref{RepForm} exists for any $z\in\Pi^2$.  Write $z_1=x_1+iy_1,\ z_2=x_2+iy_2$, with $x_1,x_2\in\R$ and $y_1,y_2 > 0$ and let $T=B + z_1Y + z_2(1-Y)$.  We have, for any $u\in\M$,
\begin{align*}
\im\ip{Tu}{u} &= y_1\ip{Yu}{u}+y_2\ip{(1-Y)u}{u} \\
	&\geq \min\{y_1,y_2\} \|u\|^2,
\end{align*}
and therefore
\[
\|Tu\| \, \|u\| \geq |\ip{Tu}{u}| \geq \im\ip{Tu}{u}\geq \min\{y_1,y_2\} \ \|u\|^2.
\]
Thus the operator $T$ has the positive lower bound $\min\{y_1,y_2\}$, and so has a left inverse.  A similar argument with $z_j$ replaced by its complex conjugate shows that $T^*$ also has a left inverse.  Hence $B + z_1Y + z_2(1-Y)$ is invertible for all $z\in\Pi^2$, and clearly the two-variable resolvent $(B + z_1Y + z_2(1-Y))^{-1}$ is analytic on $\Pi^2$.

(2)$\Rightarrow$(1)   Suppose that a representation of the form \eqref{RepForm} holds for $h$.  Then $h$ is analytic on $\Pi^2$.  For any invertible operator $T$, $\im (T^{-1})$ is congruent to $-\im T$, and so
\[
\im  (B + z_1Y + z_2(1-Y))^{-1} \mbox{ is congruent to } -(\im z_1)Y-(\im z_2)(1-Y).
\]
Since the last operator is negative, it follows that $\im h(z) \geq 0$ for all $z\in\Pi^2$, and so $h\in\p$.

To see that $\infty$ is a carapoint for $h$ note that
\[
y\im h(iy\chi) = -y\im \ip{(B+iy)^{-1}\al}{\al}.
\]
Now
\[
\im(B+iy)^{-1} = -y (B+iy)^{-1}(B-iy)^{-1}.
\]
Let the spectral representation of $B$ be 
\[
B=\int t \ \dd E(t).
\]
Then
\begin{align*}
y\im h(iy\chi) &=  y^2 \ip{(B+iy)^{-1}(B-iy)^{-1}\al}{\al} \\
	&=y^2 \int\frac{1}{(t+iy)(t-iy)}\ip{\dd E(t)\al}{\al} \\
	&= \int \frac{y^2}{t^2+y^2} \ip{\dd E(t)\al}{\al} \\
	&\to \int \ip{\dd E(t)\al}{\al} = \|\al\|^2 \quad \mbox{ as } y\to \infty
\end{align*}
by the Dominated Convergence Theorem.  Hence $\infty$ is a carapoint for $h$ with finite value.

(1)$\Rightarrow$(2)  Suppose that (1) holds and let $\ph\in\s$ be the Schur-class function associated with $h$ by equations \eqref{phiandh}.  By Proposition \ref{carainfty}, $\chi = (1,1)$ is a carapoint for $\ph$ and $\ph(\chi)\neq 1$.

By Theorem \ref{mainExistence} there exists a generalized model $(\M, u,I)$ of $\vp$  having $\chi$ as a $C$-point and an accompanying unitary realization $(a, \beta, \gamma, Q)$ of $(\M,u,I)$  with $\ker (1-Q) = \{0\}$. 
Moreover $I$ is expressible by the formula \eqref{defI} (with $\tau_1=\tau_2=1$) for some positive contraction $Y$ on $\M$.
 Thus
\[
 L = \begin{bmatrix}
      a & 1 \otimes \beta \\ \gamma \otimes 1 & Q
     \end{bmatrix}
\]
is unitary on $\C\oplus\M$ and 
\beq\label{propL}
 L\vectwo{1}{I(\lambda)u_\la} = \vectwo{\vp(\lambda)}{u_\la}.
\eeq
We wish to define the Cayley transform  $J$ of $L$:
\[
 J = i\frac{1 + L}{1 - L}.
\]
Of course $1 - L$ may not be invertible, and so we define $J$ as an operator from $\ran(1-L)$ to $\ran(1+L)$ by
\beq\label{defJ}
 J (1-L)x = i(1+L)x.
\eeq
This equation does define  $J$ as an operator, in view of the following observation.
\begin{proposition}
 If $\chi$ is a $B$-point for $\vp$ such that $\vp(\chi) \neq 1$ and $(a, \beta, \ga,Q)$ is a realization of a generalized model of $\ph$ such that $\ker (1 - Q) = \{0\}$, then $\ker(1 -L) = \{0\}$. 
\end{proposition}
\begin{proof}
 Let $x\in \ker(1 -L)\subset \C\oplus \M$ and suppose $x\neq 0$.   Since $\ker(1-Q) = \{0\}$, it cannot be that $x\in\M$, and so we can suppose that $x=\vectwo{1}{x_0}$ for some $x_0\in\M$. Then
\[
 \begin{bmatrix}
  a & 1 \otimes \beta \\ \gamma \otimes 1 & Q
 \end{bmatrix} \vectwo{1}{x_0} = \vectwo{1}{x_0},
\]
which implies that
\begin{align}\label{axb1}
  a + \ip{x_0}{b} &= 1 \\
  \gamma + Qx_0 &= x_0, \nn
\end{align}
and hence
\beq\label{x0ga}
 (1 - Q)x_0 = \gamma.
\eeq
By equation \eqref{propL}, 
\beq\label{uQI}
 u_\la = \gamma + QI(\lambda)u_\la.
\eeq
 Since $\chi$ is a $C$-point of the generalized model $(\M,u,I)$, there is a vector $u_\chi\in\M$ such that $u_\la\to u_\chi$ as $\la\nt\chi$.  On taking nontangential limits in equation \eqref{uQI} we obtain
\[
 u_\chi = \gamma + Qu_\chi,
\]
and so 
\beq\label{chiga}
(1 - Q)u_\chi = \ga.
\eeq
On comparing this relation with equation \eqref{x0ga} and using the fact that $\ker(1 - Q)=\{0\}$ we deduce that $x_0 = u_\chi$.  Again by equation \eqref{propL},
\[
  \vp(\lambda) = a + \ip{I(\lambda)u_\la}{\beta}.
\]
Let $\la\nt\chi$: then $I(\la)\to 1$ and so
\[
 \vp(\chi) = a + \ip{u_\chi}{\beta} = a +\ip{x_0}{\beta}.
\]
In view of equation \eqref{axb1} we have $\vp(\chi) = 1$, contrary to hypothesis. Thus $\ker (1 - L) = \{0\}$.
\end{proof}
We have shown that $J:\ran(1-L)\to \C\oplus\M$ is well defined by equation \eqref{defJ}. Moreover $\ran (1 - L)$ is dense in $\C\oplus\M$, since 
\[
\ran(1-L)^\perp = \ker(1-L\ad) = \ker(1-L) = \{0\}.
\] 
Thus $J$ is a densely defined operator on $\C\oplus\M$,  and since $L$ is unitary,  $J$ is self-adjoint.

The next step is to derive a matricial representation of $J$ on $\C\oplus\M$. By the definition \eqref{defJ} of $J$ and equation \eqref{propL},
\begin{align*}
 J(1-L)\vectwo{1}{I(\lambda)u_\lambda} &= i(1+L)\vectwo{1}{I(\lambda)u_\lambda}
\end{align*}
and therefore
\begin{align*}
 J\vectwo{1 - \vp(\lambda)}{(I(\lambda) - 1)u_\lambda} &= i\vectwo{1 + \vp(\lambda)}{(I(\lambda) + 1)u_\lambda}.
\end{align*}
Divide through by $1 - \vp(\lambda)$ to get
\beq\label{propJ} \
 J\vectwo{1}{\ds \frac{I(\lambda) - 1}{1 - \vp(\lambda)}u_\lambda} = \vectwo{\ds i\frac{1 + \vp(\lambda)}{1 - \vp(\lambda)}}{\ds i \frac{I(\lambda)+1}{1 - \vp(\lambda)}u_\lambda}. 
\eeq
Define $v:\Pi^2 \to \M$ by
\beq\label{defv}
v_z = - \frac{I(\lambda)-1}{1 - \vp(\lambda)}u_\lambda. 
\eeq
Recall that (compare equation \eqref{Iminus1})
\[
I(\la)-1 = -\frac{(\la_1-1)(\la_2-1)}{1-\la_1(1-Y)-\la_2 Y},
\]
and hence $I(\la)-1$ is invertible for $\la\in\D^2$.  We have
\begin{align*}
 i \frac{I(\lambda)+1}{1 - \vp(\lambda)}u_\lambda &= i \frac{I(\lambda) + 1}{I(\lambda) - 1} \left[ \frac{I(\lambda) - 1}{1 - \vp(\lambda)}u_\lambda \right] \\
 &= i \frac{1 + I(\lambda)}{1 - I(\lambda)} v_z.
\end{align*}
A straightforward calculation now yields the appealing formula
\[
 i\frac{1 + I(\lambda)}{1 - I(\lambda)} = z_1Y +z_2(1-Y).
\]
Thus equation \eqref{propJ} becomes
\beq\label{goodJ}
 J\vectwo{1}{-v_z} = \vectwo{h(z)}{\left(z_1Y +z_2(1-Y)\right)v_z}.
\eeq
We wish to write $J$ as an operator matrix
\beq\label{Jma}
 J = \begin{bmatrix}
      b & 1 \otimes \alpha \\ \alpha \otimes 1 & B
     \end{bmatrix}
\eeq
on $\C\oplus\M$, but in order for this to make sense we require that $\vectwo{1}{0}$ be in the domain of $J$, which is $\ran (1-L)$.  We must show that there exists a vector $\vectwo{c}{x}$ such that 
\[
 \begin{bmatrix}
  1 - a & -1 \otimes \beta \\ -\gamma \otimes 1 & 1 - Q
 \end{bmatrix} \vectwo{c}{x} = \vectwo{1}{0},
\]
which is to say that there exist $c\in\C$ and $x\in\M$ such that
\begin{align}\label{cx}
  c(1-a) - \ip{x}{\beta} &= 1, \\
 -c\gamma + (1-Q)x &= 0. \nn
\end{align}
Since $\ph(\chi)\neq 1$ we may choose
\[
c = \frac{1}{1 - \vp(\chi)}, \qquad x = cu_\chi,
\]
and by virtue of equation \eqref{chiga},  $c, \ x$ then satisfy equations \eqref{cx}.  Accordingly equation \eqref{Jma} is a {\em bona fide} matricial representation of $J$ on $\C\oplus\M$ for some $b\in\R$, some $\al\in\M$ and some operator $B$ on $\M$.  One can show that in fact $B$ is a densely defined self-adjoint operator on $\M$; the details can be found in, for example, \cite[Lemma 6.24]{AMY11}.

Equation \eqref{goodJ} now becomes
\[
  \begin{bmatrix}
      b & 1 \otimes \alpha \\ \alpha \otimes 1 & B
     \end{bmatrix} \vectwo{1}{-v_z} = J\vectwo{1}{-v_z} = \vectwo{h(z)}{\left(z_1Y +z_2(1-Y)\right)v_z}
\]
and so
\begin{align*}
 h(z) &= b - \ip{v_z}{\alpha}, \\
 \left(z_1Y + z_2(1 - Y)\right)v_z &= \alpha - Bv_z.
\end{align*}
Thus
\[
 v_z = \left(B + z_1Y + z_2(1-Y)\right)^{-1}\alpha,
\]
and finally 
\[
 h(z) = b - \ip{\left(B + z_1Y + z_2(1-Y)\right)^{-1}\alpha}{\alpha}.
\]
Therefore (1)$\Rightarrow$(2).
\end{proof}
Some generalizations of Nevanlinna's representation theorems to several variables can be found in \cite{ATY}.

\nin J. Agler, Department of Mathematics, University of California at San Diego, CA 92103, USA.\\

\nin R. Tully-Doyle,  Department of Mathematics, University of California at San Diego, CA 92103, USA.\\

\nin N. J. Young, School of Mathematics, Leeds University, Leeds LS2 9JT {\em and} School of Mathematics and Statistics, Newcastle University, Newcastle upon Tyne NE3 4LR, England.  Email N.J.Young@leeds.ac.uk

\end{document}